\documentclass[11pt]{article}
\usepackage{graphicx, color}
\usepackage{graphicx}

\def\la{\bigl\langle}
\def\ra{\bigr\rangle}

\def\ds{\displaystyle}
\def\forall{\hbox{for all}~}
\def\L{{\bf L}}

\def\avint{-\!\!\!\!\!\!\int}

\def\bfn{{\bf n}}

\def\ve{\varepsilon}

\def\wto{\rightharpoonup}

\def\A{{\cal A}}
\def\E{{\cal E}}
\def\I{{\cal I}}

\def\S{{\cal S}}
\def\H{{\cal H}}

\def\R{I\!\!R}

\def\implies{\Longrightarrow}
\def\vp{\varphi}

\def\vs{\vskip 2em}

\def\v{\vskip 1em}

\def\begi{\begin{itemize}}
\def\endi{\end{itemize}}
\def\C{{\cal C}}
\def\cL{{\cal L}}

\def\M{{\cal M}}
\def\ov{\overline}
\def\Tilde{\widetilde}
\def\Hat{\widehat}
\def\bega{\begin{array}}
\def\enda{\end{array}}
\def\meas{\hbox{meas}}

\def\bel{\begin{equation}\label}
\def\eeq{\end{equation}}
\def\sqr#1#2{\vbox{\hrule height .#2pt
\hbox{\vrule width .#2pt height #1pt \kern #1pt
\vrule width .#2pt}\hrule height .#2pt }}
\def\square{\sqr74}
\def\endproof{\hphantom{MM}\hfill\llap{$\square$}\goodbreak}

\newtheorem{thm}{Theorem}[section]

\newtheorem{lemma}{Lemma}[section]

\newtheorem{remark}{Remark}[section]
\newtheorem{definition}{Definition}[section]

\begin{document}
\title{\bf  On the Optimal Shape of Tree Roots and Branches}
\vs

\author{Alberto Bressan and Qing Sun\\
\,
\\
Department of Mathematics, Penn State University, \\
e-mails: axb62@psu.edu, ~qxs15@psu.edu}
\maketitle
\begin{abstract} This paper introduces two classes of variational problems, determining
optimal shapes for tree roots and branches.
Given a measure $\mu$, describing the distribution of leaves, we introduce
a  sunlight functional $\S(\mu)$ computing the total 
amount of light captured by the leaves.
On the other hand, given a measure $\mu$ describing the distribution of root hair cells,
we consider a harvest functional $\H(\mu)$ computing the total amount of water and nutrients
gathered by the roots.  In both cases, we seek to maximize these functionals subject to a 
ramified transportation cost, for transporting 
nutrients from the roots to the trunk and from the trunk to the leaves. 
The main results establish various properties of these functionals, and the existence of
optimal distributions. In particular, we prove the upper semicontinuity 
of $\S$ and $\H$,  together with a priori estimates on the  support of optimal distributions.

\end{abstract} 
\vs
\section{Introduction}
\label{s:0}
\setcounter{equation}{0}
Living organisms come in 
an immense variety of shapes, such as
roots, branches, leaves, and flowers in plants, or bones in animals.  
In many  cases, it is expected that through natural selection, 
these organisms have evolved into a ``best possible" shape.
From a mathematical perspective, it is thus  of interest to 
study functionals whose minimizers may determine
some of the many shapes found in the biological world.

As a step in this direction, 
in this paper we consider two functionals, defined on 
a space of positive measures on $\R^d$, 
and show how they can be used
to describe the optimal configurations of roots and branches in a tree.

The first one, which we call the ``sunlight functional", models  the total amount of
sunlight captured by the leaves of  a tree.
  Here we think of  a measure $\mu$ 
as the density of leaves.   To achieve a realistic model, 
 our functional $\S(\mu)$ will take different forms in the case of a free-standing 
 tree in the middle of a prairie, or a tree in a forest, whose lower
branches are partially shielded by the surrounding vegetation. 
The model also accounts for the fact that 
light rays come from different directions at different times of the day.
 
The second one, which we call the ``harvest functional", models the 
total amount of water and nutrients collected by the roots.   
In this case,  we think of a measure $\mu$ as the density of root  hair cells in the soil.
A similar harvest functional $\H(\mu)$  was introduced in \cite{BCS},
in connection with a problem of optimal harvesting of marine resources.
In the present paper, both Dirichlet and Neumann boundary conditions will be considered.

The above functionals will be combined with a  ``ramified transportation cost", 
for transporting nutrients from the roots to the base of the trunk, 
or from the base of the trunk to the leaves.  For a given measure $\mu$
on $\R^d$, this is modeled by the minimum $\alpha$-irrigation cost 
$\I^\alpha(\mu)$ from the 
origin, introduced in \cite{MMS, X3}.  The lower semicontinuity of this cost
plays an essential role toward the existence of optimal solutions.
For a comprehensive introduction to 
optimal irrigation problems  we refer to \cite{BCM}.

The optimal shape of branches is now determined by the variational problem
\bel{msun1}\hbox{maximize:}\quad \S(\mu)-c\I^\alpha(\mu)\eeq
for some constants $0<\alpha<1$ and $c>0$.  We  study this maximization problem
among all positive measures with a given total mass: 
\bel{tmass}\mu(\R^d)\, = \,\kappa_0\,.\eeq
Notice that, to maximize the gathered sunlight, the leaves should  be  spread
out as wide as possible.  On the other hand, this makes it more costly to transport nutrients
from the root to all the leaves.   

Similarly, the optimal structure of a root system  can be related  to the problem
\bel{mh1}\hbox{maximize:}\quad \H(\mu)-c\I^\alpha(\mu).\eeq

The remainder of the paper is organized as follows.   In Section~2 we introduce
a sunlight functional and prove some of its  properties.  These include
the upper semicontinuity and various estimates.  Section~3 is concerned with 
the harvest functional, recalling the main definitions and extending
some of the results in \cite{BCS} to different boundary conditions. 
In Section~4 we briefly review the theory of optimal ramified transport, proving some
estimates on
the minimum $\alpha$-irrigation cost for a measure $\mu$, for later use.
The optimization problems for the shape of tree branches
and tree roots are studied in Sections~5 and 6, respectively. 
Using the semicontinuity of the various functionals, together with a priori 
bounds on the supports of a sequence of optimizing measures, in both cases
we establish the existence of an optimal solution.  Some concluding remarks
are given in the last section.
\v
\section{The sunlight functional}
\label{s:2}
\setcounter{equation}{0}

Throughout the following, $B(x_0, r) $ denotes an open ball centered at 
$x_0$ with radius $r$, while
$S^{d-1}=\{x\in\R^d\,;~|x|=1\}$ denotes the unit sphere in $\R^d$.
We write $\ov \Omega$ for the closure of a set $\Omega$, and  $\cL^d$ for the 
$d$-dimensional Lebesgue measure.

Let $\mu$ be a positive, bounded Radon measure on $\R^d$.   
Thinking of $\mu$ as the distribution of leaves on a tree, 
we seek a functional $\S(\mu)$ describing 
the total amount of sunlight captured by the leaves.

\begin{figure}[ht]
\centerline{\hbox{\includegraphics[width=10cm]{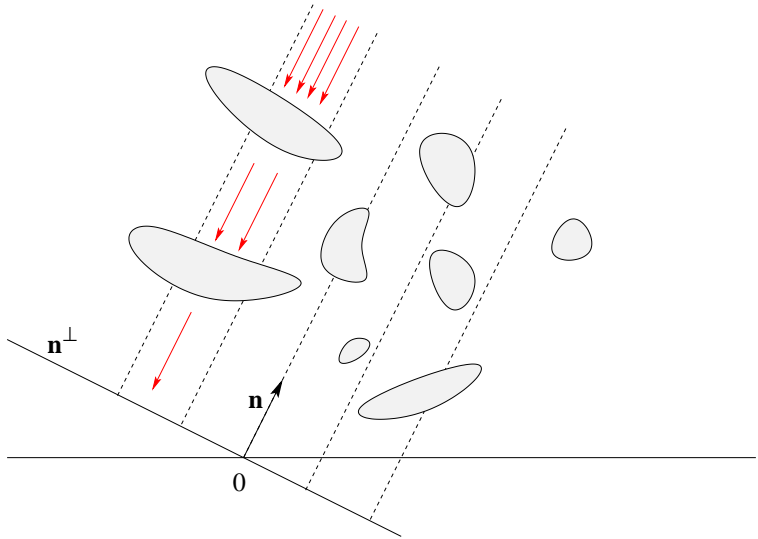}}}
\caption{\small  Sunlight arrives from the direction 
parallel to $\bfn$.  
Part of it is absorbed by the measure $\mu$, 
supported on the grey regions.}
\label{f:ir10}
\end{figure}

To begin with a simple setting, fix a unit vector
$\bfn\in \R^d$ 
and assume that all  
light rays come parallel to $\bfn$.  Moreover, assume 
that the measure
$\mu$ is absolutely continuous with density 
$f$ w.r.t.~Lebesgue measure on $\R^d$.  
Call $E_\bfn^\perp$ the $(d-1)$-dimensional subspace perpendicular to $\bfn$ 
and let $\pi_\bfn:\R^d\mapsto \bfn^\perp$ be the perpendicular projection.   As shown in 
Fig.~\ref{f:ir10}, each point $ x\in \R^d$ can  be expressed 
uniquely as
\bel{perp}
 x~=~ y + s\bfn\eeq
with $ y\in E_\bfn^\perp$ and $s\in\R$.

Our basic modeling assumption is that the rate at which sunlight is absorbed is proportional
to the local density of leaves.   
For each fixed $ y\in E^\perp_\bfn$, calling
$s\mapsto \phi( y,s)$ 
the amount of sunlight reaching the point 
$ x =  y + s\bfn$, we thus assume
$${\partial\over\partial s}\phi( y,s)~=~f( y + s\bfn) \phi( y,s)$$
$$\lim_{s\to +\infty} \phi( y,s)~=~1.$$
For simplicity, we here assign unit values to the absorption rate, and to the amount of light 
arriving from infinity per unit  $(d-1)$-dimensional volume in $E^\perp_\bfn$.
This implies
\bel{phi4}\phi( y,s)~=~
\exp\left\{ -\int_s^{+\infty}f( y + t\bfn)\, dt \right\}.\eeq
Integrating over the perpendicular plane $E_\bfn^\perp$, the total amount of light which is absorbed by the leaves
is thus
\bel{sun1}
\S^\bfn(\mu)~=~\int_{E_\bfn^\perp}~1~-~\exp\left\{ -
\int_{-\infty}^{+\infty}f( y + t\bfn)\, dt \right\}\, 
d y.\eeq
We now observe that the formula (\ref{sun1}) can be easily 
extended to the case of
a general measure $\mu$, not necessarily absolutely continuous w.r.t.~Lebesgue measure.

On the perpendicular subspace $E_\bfn^\perp$ consider the projected 
measure $\mu^\bfn$ defined by setting
\bel{mupro}\mu^\bfn(A)~=~\mu\Big(\bigl\{ x\in\R^d\,;~~\pi_\bfn(x)\in A\bigr\}\Big)\eeq
for every open set $A\subseteq E_\bfn^\perp$.
Call $\Phi^\bfn$ the density of the absolutely continuous part of $\mu^\bfn$
w.r.t.~the $(d-1)$-dimensional Lebesgue measure on $E_\bfn^\perp$.
\v
\begin{definition}
 The total amount of sunlight from the direction $\bfn$ absorbed by a measure 
$\mu$ on $\R^d$ is defined as
\bel{SSn}
\S^\bfn(\mu)~\doteq~
\int_{E_\bfn^\perp}\Big(1- \exp\bigl\{ - \Phi^\bfn(y)\bigr\}\Big)
\, dy\,.\eeq
\end{definition}

Next, we model the fact that sunlight does not always come from the same direction.  Instead, there exists a density function $\eta:S^{d-1}\mapsto \R_+$
which describes the total amount of light coming from the direction $\bfn$
during the course of a day.
\v
\begin{definition} If light comes from different directions with variable intensity 
$\eta=\eta(\bfn)$, 
the total amount of sunshine  captured by a measure 
$\mu$ on $\R^d$ is defined as
\bel{SSe}
\S^\eta(\mu)~\doteq~
\int_{S^{d-1}} \S^\bfn(\mu)\,\eta(\bfn)\,d\bfn\,.\eeq
\end{definition}
\v
\begin{remark}  {\rm Measures
which are singular w.r.t.~the $(d-1)$-dimensional Hausdorff measure
are irrelevant.  More precisely, if  $\mu= \mu_1+\mu_2$ and $\mu_2$ is supported on a set whose $(d-1)$-dimensional measure is zero, then
$\S^\bfn(\mu) = \S^\bfn(\mu_1)$ for every unit vector $\bfn\in\R^d$.
}\end{remark}

\begin{remark}  {\rm A case of particular interest is when light comes uniformly
from all directions of the positive half sphere.
$$S^{d-1}_+~\doteq~\Big\{ \bfn  =(n_1, \ldots, n_d)\,;~~|\bfn|=1\,,~~n_d>0\Big\},$$
We shall model this situation  by taking 
\bel{eta1}\eta(\bfn) ~=~\left\{\bega{cl}  \sigma_d/2\quad\hbox{if}\quad
\bfn\in S^{d-1}_+\,,\\[4mm] 
  0\quad\hbox{otherwise.}\enda
\right.\eeq
Here $\sigma_d$ denotes the $(d-1)$-dimensional measure of the 
surface of the unit ball in $\R^d$.}
\end{remark}
\v
The next lemma collects some elementary properties of the functional $\S^\bfn$.
In the following, we denote by $\mu^\lambda$ the measure such that
\bel{mldef}\mu^\lambda(A)~=~\mu(\lambda^{-1} A)\eeq
for  every open set $A\subset \R^d$, so that
\bel{suppml}
Supp(\mu^\lambda)~=~\lambda\cdot Supp(\mu)~=~\{\lambda x\,;~~x\in Supp(\mu)\}.\eeq
Moreover, $\omega_d$ denotes the volume of the unit ball in $\R^d$.
\v
\begin{lemma}\label{l:sl} Let $\mu,\tilde \mu$ be positive Radon measures on $\R^d$.
For any unit vector $\bfn\in S^{d-1}$, the following  holds.
\begi
\item[(i)] $\S^\bfn(\mu)~\leq~ \mu(\R^d)$,

\item[(ii)] If the measure $\mu$ is supported 
inside a ball of radius $r$, then 
$
\S^\bfn(\mu)\leq  \omega_{d-1}\, r^{d-1}$.

\item[(iii)] 
$\S^\bfn(\mu)~\leq~\S^\bfn(\mu+\tilde\mu)~\leq~\S^\bfn(\mu)+\S^\bfn(\tilde\mu)$.

\item[(iv)] $\S^\bfn(\lambda\mu)~\leq~\lambda S^\bfn(\mu)$, 
for every $\lambda\geq 1$.

\item[(v)] For every $\lambda>0$ one has
\bel{Sdial}\S^\bfn(\lambda^{d-1}\mu^\lambda) ~=~\lambda^{d-1}\S^\bfn(\mu).\eeq

\item[(vi)] If $\mu$ is absolutely continuous w.r.t.~Lebesgue measure, then
\bel{lre}
\lim_{\lambda\to 0+}  {\S^\bfn(\lambda\mu)\over\lambda}~
=~\mu(\R^d)~=~\lim_{\lambda\to +\infty} \S^\bfn(\mu^\lambda).\eeq
\endi
\end{lemma}
\v
{\bf Proof.} 
{\bf 1.} To prove (i), consider any unit vector $\bfn$
and call $\Phi^\bfn$ the density of the absolutely continuous part of of the projected measure $\mu^\bfn$
w.r.t.~the $(d-1)$-dimensional Lebesgue measure on $E_\bfn^\perp$.
Then 
\bel{Snbo}\S^{\bfn}(\mu)~\doteq ~\int_{E_\bfn^{\perp}}\Big( 1-\exp\{-\Phi^{\bfn}(y)\}\Big)\,dy~\leq~ 
\int_{E_\bfn^{\perp}}\Phi^{\bfn}(y)\,dy~\leq ~\mu^\bfn(E^\perp_\bfn)~=~\mu(\R^{n}).\eeq
\v
{\bf 2.} To prove (ii), let $\mu$ be supported inside the ball $B(x_0, r)$, centered at $x_0$ with radius $r$.
Call $y_0$ the perpendicular projection of $x_0$ on the space $E^\perp_\bfn$.
Then $\Phi^\bfn(y)=0$ whenever $|y-y_0|>r$. Hence
$$\S^{\bfn}(\mu)~\doteq ~\int_{E_\bfn^{\perp}}\Big( 1-\exp\{-\Phi^{\bfn}(y)\}\Big)\,dy~\leq~ 
\int_{E_\bfn^{\perp}\cap B(y_0,r)} 1\, dy~=~\omega_{d-1}r^{d-1}.$$
\v
{\bf 3.} Concerning (iii), the first inequality is an immediate consequence of 
the monotonicity of the function $1-e^{-x}.$ 
	 To prove the second inequality, denote by  $\Phi^{\bfn},\Tilde \Phi^{\bfn} $ the density functions of  the projected measures $\mu^\bfn, \tilde \mu^\bfn$ on the perpendicular space
	 $E^\perp_\bfn$.  Observing that $\Phi^{\bfn}+\Tilde \Phi^{\bfn} $  
	 is the density function of $(\mu+\tilde \mu)^\bfn$, one obtains
	\bel{lss3}
	\begin{array}{l}
	\ds
	 \S^{\bfn}(\mu+\tilde \mu)-\S^{\bfn}(\mu)-\S^{\bfn}(\tilde \mu)\\[4mm]
	\qquad \ds=~\int_{E_\bfn^\perp}\Big(1-\exp\bigl\{-\Phi^{\bfn}(y)-\Tilde \Phi^{\bfn}(y)
	\bigr\} \Big)-\Big(1-\exp\{-\Phi^{\bfn}(y)\}\Big)-\Big(1-\exp\{-\Tilde\Phi^{\bfn}(y)\}\Big)\,dy 
	\\[4mm]
	\ds\qquad=~\int_{E_\bfn^{\perp}}\Big[\exp\{-\Phi^{\bfn}(y)\}+\exp\{-\Tilde\Phi^{\bfn}(y)\}
	-\exp\{-\Phi^{\bfn}(y)-\Tilde\Phi^{\bfn}(y)\}-1\Big]\,dy~\leq~0.
	\end{array}
	\eeq
	Indeed, the last inequality is obtained by checking that 
	$$h(x_1,x_2)~\doteq ~e^{-x_1}+e^{-x_2}-e^{-x_1-x_2}-1~\leq~ 0$$
	for every $x_1,x_2\geq 0$.

{\bf 4.} To prove (iv),
consider the function 
 $$h(x)~\doteq ~1-e^{-\lambda x}-\lambda +\lambda e^{- x}.$$
Assuming $\lambda\geq 1$, an elementary computation yields
$$h(0)~=~0,\qquad\qquad 
h'(x)~=~\lambda\, e^{-\lambda x}-\lambda\, e^{-x}~\leq~ 0\qquad \hbox{for all}~~x\geq 0.$$
Therefore $h(x)\leq 0$ for all $x\geq 0$.   Using this inequality we obtain
    \bel{lss2}
    \begin{array}{rl}
    S^{\bfn}(\lambda\mu)-\lambda S^{\bfn}(\mu)&=~\ds\int_{E_\bfn^\perp}\Big( 1- \exp\{-\lambda\Phi^{\bfn}(y)\}\Big)dy-\int_{E_\bfn^\perp}\lambda\Big(1-\exp\{-\Phi^{\bfn}(y)\}\Big)dy\\[4mm]
   &=~\ds\int_{E_\bfn^\perp} h\bigl(\Phi^{\bfn}(y)\bigr) \,dy~\leq~0.
   \end{array}\eeq
\v	

\v
{\bf 5.} 
To prove (v),  we first compute the density function $\Phi^{\bfn,\lambda}$ for
the projected  measure $(\lambda^{d-1}\mu^\lambda)^\bfn$ on the $(d-1)$-dimensional subspace
$E^\perp_\bfn$.   From the identity
$$\int_A\Phi^{\bfn,\lambda}(y)\, dy~=~\int_{\lambda^{-1} A} \lambda^{d-1}\Phi^{\bfn}(\tilde y)\, d\tilde y
~=~ \int_A  \Phi^{\bfn}(\lambda^{-1} y)\, dy
$$
valid for every open set $A\subseteq E^\perp_\bfn$, we deduce
$$\Phi^{\bfn,\lambda}(y)~=~\Phi^\bfn(\lambda^{-1} y)$$
for every $y\in E^\perp_\bfn$.
    Therefore, using the change of variable  $\tilde y = \lambda^{-1} y$, one obtains 
    \begin{equation}
    \begin{array}{rl}
    \S^\bfn(\lambda^{d-1}\mu^\lambda)&\ds=~\int_{E_\bfn^\perp}\Big(1-\exp\{-\Phi^{\bfn,\lambda}(y)\}\Big)
    \,dy\\[4mm]
 	&=~\ds\int_{E_\bfn^\perp}\lambda^{d-1}\Big(1-\exp\{-\Phi^\bfn(\tilde y)\}\Big)\,d\tilde y\\[4mm]
 	&=~\ds\lambda^{d-1}\S^\bfn(\mu).    \end{array}
    \end{equation}
\v
{\bf 6.}   It remains to prove the
two limits in (\ref{lre}).  Assume that the  positive measure $\mu$ has density
$f$ w.r.t.~Lebesgue measure on $\R^d$. Then
$$\frac{\S^\bfn(\lambda\mu)}{\lambda}~=~\int_{E^{\perp}_{\bfn}}\frac{\Big(1-\exp\Big\{-\lambda\int_{-\infty}^{\infty}f(y+t\bfn)dt \Big\}\Big)}{\lambda}~dy$$
By Fubini's theorem, for almost every $y\in E_{\bfn}^\perp$ we have $\int_{-\infty}^\infty f(y+t\bfn)~dt<\infty$.
 At such a point $y$ we have
$$\lim_{\lambda\to 0+}\frac{1-\exp\Big\{-\lambda\int_{-\infty}^\infty f(y+t\bfn)\Big\}}{\lambda}~=~\int_{-\infty}^{\infty}f(y+t\bfn)~dt.$$
On the other hand,
\bel{dbd1}
\frac{1-\exp\Big\{-\lambda\int_{-\infty}^\infty f(y+t\bfn)\Big\}}{\lambda}~\leq~\int_{-\infty}^{\infty}f(y+t\bfn)~dt
\eeq
Therefore, by dominated convergence theorem, we conclude
\bel{dbd2}
\lim_{\lambda\to 0+}\frac{\S^\bfn(\lambda\mu)}{\lambda}~=~\int_{E^\perp_\bfn}\int_{-\infty}^\infty f(y+t\bfn)~dtdy~=~\mu(\R^d).
\eeq

To prove the second equality in (\ref{lre}), call $\Phi^{\bfn,\lambda}$ the density function 
for the projected measure $(\mu^\lambda)^\bfn$.  For almost every $y\in E^\perp_\bfn$ we have
\bel{Pnn}
\Phi^{\bfn,\lambda}(y)~=~\frac{1}{\lambda^{d-1}}\int_{-\infty}^\infty f\left(\frac{y}{\lambda}+t\bfn\right) dt~<~+\infty.\eeq
Therefore, 
\bel{wrpenn}
\bega{l}
\ds~\S^\bfn(\mu^\lambda)~=~\int_{E^\perp} \left(1-\exp\Big\{-\frac{1}{\lambda^{d-1}}\int_{-\infty}^\infty f\Big(
\frac{y}{\lambda}+t\bfn\Big)dt\Big\}\right)\,dy\\[4mm]
\ds~\quad\qquad~=~\int_{E^\perp}\lambda^{d-1}\Big(1-\exp\Big\{-\frac{1}{\lambda^{d-1}}\int_{-\infty}^\infty f(y+t\bfn)dt\Big\}\Big)~dy.
\enda
\eeq
For a.e.~$y$ we have $\int_{-\infty}^\infty f(y+t\bfn)dt<\infty$, and hence
$$\lambda^{d-1}\left(1-\exp\left\{-\frac{1}{\lambda^{d-1}}\int_{-\infty}^\infty f(y+t\bfn)dt\right\}\right)~\leq~\int_{-\infty}^\infty f(y+t\bfn)\,dt.$$
On the other hand, by L'Hospital Rule, 
$$\lim_{\lambda\to +\infty}\lambda^{d-1}\Big(1-\exp\Big\{-\frac{1}{\lambda^{d-1}}\int_{-\infty}^\infty f(y+t\bfn)dt\Big\}\Big)~=~\int_{-\infty}^\infty f(y+t\bfn)~dt$$
Letting $\lambda\to +\infty$ in  (\ref{wrpenn}), by 
the dominated convergence theorem one obtains the second  equality in (\ref{lre}).
\endproof
	


\v
The formula (\ref{SSn})  covers the case where 
there are no other obstacles to light propagation except $\mu$.
Next, we want to model the presence of other plants that 
capture part of the light, and determine how much light 
is actually collected by $\mu$.

As a preliminary, consider two positive  measures 
$\mu$ and $\nu$,
absolutely continuous  with densities $f,g$ 
w.r.t.~Lebesgue measure on $\R^d$.   Assuming that 
light comes from the direction  $\bfn$, the same 
computation as in (\ref{phi4})
shows that the total 
amount of light that reaches a 
point $ x =  y + s\bfn$ is
$$\exp\left\{-\int_s^{+\infty} \bigl( f( y+t\bfn) + g( y+t\bfn)
\bigr)\,dt\right\}.$$
Integrating by parts, the total amount of light collected by 
the distribution $\mu$ with density $f$ is computed by
\bel{If1}\bega{rl}
\S^\bfn(\mu)&\ds=~\int_{E_\bfn^\perp}\left(
\int_{-\infty}^{+\infty} f( y+s\bfn)\,\exp\left\{ -\int_s^{+\infty} \bigl(f( y+t\bfn)+g( y+t\bfn)\bigr)\, dt\right\}\, ds\right)d y\\[4mm]
&\ds =~\int_{E_\bfn^\perp}\left(\int_{-\infty}^{+\infty}{d\over ds} \exp\left\{ -\int_s^{+\infty} f( y+t\bfn)\, dt\right\}\cdot
\exp\left\{ -\int_s^{+\infty} g( y+t\bfn)\, dt\right\}\, ds\right)dy\\[4mm]
&\ds =~\int_{E_\bfn^\perp}\Bigg( 
1- \exp\left\{ -\int_{-\infty}^{+\infty} f( y+t\bfn)\, dt\right\}
\exp\left\{ -\int_{-\infty}^{+\infty} g( y+t\bfn)\, dt\right\}\\[4mm]
&\ds\qquad\qquad - \int_{-\infty}^{+\infty} g( y+s\bfn) \exp\left\{ -\int_s^{+\infty} \bigl(f( y+t\bfn)+g( y+t\bfn)\bigr)\, dt\right\}\, ds\Bigg)dy\,.
\enda
\eeq
In essence, this says that 
$$\hbox{[light collected by $\mu$]} ~=~\hbox{[light collected by $\mu+\nu$]} -\hbox{[light collected by $\nu$]}.$$ 
Notice that here the right hand side makes sense also if $\mu$ is 
is an arbitrary measure, not necessarily absolutely continuous w.r.t.~Lebesgue measure.
This fact can be used to define the total sunlight 
absorbed by any positive measure $\mu$, in the presence of a second measure  $\nu$ which is absolutely continuous
with density $g(\cdot)$
w.r.t.~Lebesgue measure on $\R^d$. 

For any given a unit vector $\bfn$, we represent
$\R^d=E_\bfn\oplus E_\bfn^\perp$,  as the sum of the orthogonal spaces containing all vectors parallel and orthogonal to $\bfn$,
respectively.   We denote by $(t,y)
\in E_\bfn\oplus E_\bfn^\perp$ the variable corresponding to this decomposition.
As before,
let $\pi_\bfn:\R^d\mapsto E_\bfn^\perp$ be the perpendicular projection, and call $\mu^\bfn$ be the projection of $\mu$
on $E_\bfn^\perp$, defined as in (\ref{mupro}).
By Theorem~2.28 in \cite{AFP} (on the disintegration of the measure $\mu$), there exists a family of 1-dimensional measures
$\mu^y$, $y\in E_\bfn^\perp$, such that the following holds.

\begi
\item[(i)]  $\mu^y(E_\bfn)=1$ for every $y\in E_\bfn^\perp$.
\item[(ii)]  The map $y\mapsto \mu^y$ is $\mu^\bfn$-measurable.
\item[(iii)] 
For every $\phi\in \L^1(\R^d)$ one has
\bel{disint}
\int_{\R^d} \phi \,d\mu~=~\int_{E_\bfn^\perp}\left(
\int_{-\infty}^{+\infty}\phi(t,y)d\mu^y(t)\right)d\mu^\bfn(y)\,.\eeq
\endi

\begin{figure}[ht]
\centerline{\hbox{\includegraphics[width=6cm]{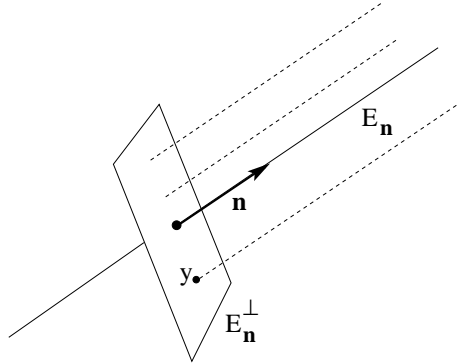}}}
\caption{\small Disintegration of a measure $\mu$ on $\R^d$.
According to (\ref{disint}), the integral $\int \phi \,d\mu$ 
can be computed first integrating $\phi$ along each line 
$\{y+t\bfn\,;\,t\in\R\}$ parallel to the unit vector 
$\bfn$, then integrating over the variable 
$y\in E_\bfn^\perp$.  }
\label{f:ir5}
\end{figure}
To compute the total amount 
of light coming from the direction parallel to $\bfn$
which is captured by the measure $\mu$, we proceed as follows.

Let $\Phi^\bfn$ be the density of the 
 absolutely continuous   part of $\mu^\bfn$
w.r.t.~$(d-1)$-dimensional Lebesgue measure on $E_\bfn^\perp$,
as in (\ref{SSn}).

Now let $\nu$ be a second measure, absolutely continuous
with density $g$ w.r.t.~Lebesgue measure on $\R^d$.
Motivated by (\ref{If1}), for each $y\in E_\bfn^\perp$
we define
\bel{Imny}\bega{rl}
\S^\bfn_{\mu,\nu}(y)&\ds\doteq~1-\exp\{-\Phi^\bfn(y)\}
\exp\left\{
-\int g(y+s\bfn)\, ds\right\} \\[4mm]&\quad \ds -
\int \left( g(y+s\bfn) \exp\left\{
-\int_s^{+\infty} g(y+t\bfn)\, dt\right\} 
\cdot \exp\Big\{ -\Phi^\bfn(y)\cdot\mu^y\bigl( [s, +\infty[ \,
\bigr)\Big\}
\right)\, ds
\,,\enda\eeq
\bel{Imn} \S^\bfn(\mu;\nu)~\doteq~\int_{E_\bfn^\perp}
\S^\bfn_{\mu,\nu}(y)\, dy\,.\eeq
\v
\begin{definition}
Assume that light comes with variable intensity $\eta(\cdot)$  from 
different directions.  The {\bf total sunlight $\S^\eta (\mu;\nu)$
absorbed by the measure $\mu$ in the presence of the
absolutely continuous measure $\nu$} is then defined as
\bel{Im}
\S^\eta (\mu;\nu)~\doteq~\int_{S^{d-1}} \S^\bfn(\mu;\nu)\,\eta(\bfn)\,
d\bfn\,.\eeq
\end{definition}

We observe that the first three estimates in Lemma~\ref{l:sl} 
remain valid in this more general situation. 
\v
\begin{lemma}\label{l:sl2} Let $\mu,\nu $ be  positive Radon measures on $\R^d$.
Assume that   $\nu$ is absolutely continuous w.r.t.~Lebesgue measure.
For any unit vector $\bfn\in S^{d-1}$, the following  holds.
\begi
\item[(i)] 
~$
\S^\bfn(\mu,\nu)~\leq~\S^\bfn(\mu)~\leq~ \mu(\R^d)$.
\item[(ii)] If  the measure $\mu$ is supported
inside a  ball of radius $r$,
then $
\S^\bfn(\mu;\nu)~\leq~\omega_{d-1}\, r^{d-1}$.
\item[(iii)] For any positive measures $\mu_1, \mu_2$ one has
\bel{qsy}
\S^\bfn(\mu_1,\nu)~\leq~\S^\bfn(\mu_1+\mu_2,\nu)~\leq~\S^\bfn(\mu_1,\nu)+\mu_2(\R^d).
\eeq\endi\end{lemma}
\v
{\bf Proof.}
{\bf 1.} Let $g$ be the density of $\nu$ w.r.t.~Lebesgue measure on $\R^d$. By
 (\ref{Imny}) we have
	\bel{Imny1}\bega{rl}
	\S^\bfn_{\mu,\nu}(y)&\ds\doteq~1-\exp\{-\Phi^\bfn(y)\}
	\exp\left\{
	-\int g(y+s\bfn)\, ds\right\} \\[4mm]&\qquad  \ds -
	\int \left( g(y+s\bfn) \exp\left\{
	-\int_s^{+\infty} g(y+t\bfn)\, dt\right\} 
	\cdot \exp\Big\{ -\Phi^\bfn(y)\cdot\mu^y\bigl( [s, +\infty[ \,
	\bigr)\Big\}
	\right)\, ds\\[4mm]
	&\ds\leq~1-\exp\{-\Phi^\bfn(y)\}
	\exp\left\{
	-\int g(y+s\bfn)\, ds\right\} \\[4mm]&\qquad  \ds -
	\int\left( {d\over ds} \exp\left\{
	-\int_s^{+\infty} g(y+t\bfn)\, dt\right\} \right)\, ds
	\cdot \exp\bigl\{ -\Phi^\bfn(y)\bigr\}\\[4mm]
	&\ds=~1-\exp\{-\Phi^\bfn(y)\}~\leq~\Phi^\bfn(y)
	\,.\enda\eeq
Integrating over $E^\perp_\bfn$ we obtain the first inequality in 
(i).    The second inequality is now a consequence of  (\ref{Snbo}).
\v
{\bf 2.} If $\mu$ is supported in a ball of radius $r$, then 
the estimate (ii) follows immediately from $\S^\bfn(\mu,\nu)\leq\S^\bfn(\mu)\leq\omega_{d-1}r^{d-1}$.
\v
	
	{\bf 3.} To prove (iii), let $\Phi^\bfn_1,\Phi^\bfn_2$, and 
	$\Phi^\bfn=\Phi_1^\bfn+\Phi_2^\bfn$  be the densities of the
	absolutely continuous parts of $\mu_1^\bfn$, $\mu^\bfn_2$, and 
	$\mu^\bfn=\mu^\bfn_1+\mu^\bfn_2$  w.r.t.~the $(d-1)$-dimensional
	Lebesgue measure on $E^\perp_\bfn$, respectively.  We claim that 
	\bel{Hin}\S^\bfn_{\mu_1+\mu_2,\nu}(y)-\S^\bfn_{\mu_1,\nu}(y)~\leq~\Phi_2^\bfn(y)\eeq
	for almost every $y\in E^\perp_\bfn$. 
Indeed,
	for a fixed $y$,  assume $\Phi_2(y)\not= 0$
	and define
	\bel{ldef}\lambda~\doteq~\frac{\Phi^\bfn_1(y)}{\Phi^\bfn_1(y)+\Phi^\bfn_2(y)}~<~1\,.\eeq
	Call $\mu_1^y, \mu_2^y$, and $\mu^y$  the probability measures on 
	the 1-dimensional space $E_\bfn$ corresponding 
	to the disintegration of $\mu_1,\mu_2$, and $\mu=\mu_1+\mu_2$, respectively. 
By (\ref{ldef}) it follows
	\bel{muy12}\mu^y~=~\lambda \mu_1^y + (1-\lambda) \mu_2^y\,.\eeq
	We now compute
	\bel{S12es}\begin{array}{l}\ds
	\S^\bfn_{\mu_1+\mu_2,\nu}(y)-\S^\bfn_{\mu_1,\nu}(y)\\[4mm]
\quad \ds	=~\Big(\exp\{-\Phi^\bfn_1(y)\}-\exp\Big\{-\Phi^\bfn_1(y)-\Phi^\bfn_2(y)\Big\}\Big)\cdot\exp\left\{-\int g(y+s\bfn)\,ds\right\}\\[4mm]
	\qquad\qquad  \ds +\int g(y+s\bfn)\cdot \exp\left\{-\int_s^\infty g(y+t\bfn)dt\right\}\\[4mm]
	\qquad\qquad\qquad \ds
	\cdot \bigg(\exp\Big\{-\Phi^\bfn_1(y)\cdot\mu^y_1[s,\infty[\Big\}-
	\exp\Big\{-\bigl(\Phi^\bfn_1(y)+\Phi^\bfn_2(y)\bigr)\,
\mu^y[s,\infty[
\Big\}\bigg)ds\\[4mm]
\quad 	\doteq~I+J.
	\end{array}\eeq
The second term in the above expression can be estimated as
 \bel{mkk4}
\begin{array}{l}
	\ds J~=~\int g(y+s\bfn)\cdot \exp\left\{-\int_s^\infty g(y+t\bfn)dt\right\}\\[4mm]
	 \ds
 	\cdot \Bigg(\exp\left\{-\frac{\lambda}{1-\lambda}\Phi^\bfn_2(y)\cdot
	\mu^y_1[s,\infty[\right\}-\exp\left\{-\frac{1}{1-\lambda}\Phi^\bfn_2(y)\cdot
\bigl(\lambda\, \mu^y_1+(1-\lambda)\, \mu^y_2\bigr)[s,\infty[
 \right\}\Bigg)ds
 \\[4mm]
~ \ds=~\int g(y+s\bfn)\cdot\exp\left\{-\int_{s}^\infty g(y+t\bfn)dt\right\}\\[4mm]
\qquad \ds
\cdot\Bigg(\exp\left \{-\frac{\lambda}{1-\lambda}\Phi^\bfn_2(y)\cdot \mu^y_1[s,\infty[\right\}
\cdot \Big( 1 - \exp\bigl\{ \Phi^\bfn_2(y)  \cdot\mu^y[s,\infty[\,\bigr\}\Big)\Bigg)ds\\[4mm]
\ds ~ \leq~\int g(y+s\bfn)\cdot\exp\Big\{-\int_s^\infty g(y+t\bfn)dt\Big\}\cdot\Phi^\bfn_2(y)~ds\\[4mm]
 ~=~\ds \left(1-\exp\left\{-\int g(y+s\bfn)ds\right\}\right)\cdot \Phi^\bfn_2(y)
\end{array}
 \eeq
Combining (\ref{S12es}) with (\ref{mkk4}) we obtain
$$
\begin{array}{l}
\S^\bfn_{\mu_1+\mu_2,\nu}(y)-\S^\bfn_{\mu_1,\nu}(y)~=~I+J\\[4mm]
\ds\leq~\Big(\exp\{-\Phi^\bfn_1(y)
\}-\exp\{-\Phi^\bfn_1(y)-\Phi^\bfn_2(y)\}-\Phi^\bfn_2(y)\Big)\cdot 
\exp\left\{-\int g(y+s\bfn)ds\right\}+\Phi^\bfn_2(y)\\[4mm]
\leq~\Phi^\bfn_2(y).
\end{array}
$$
	Integrating over the $(d-1)$-dimensional space $E_\bfn^\perp$ 
	one obtains the desired estimate. This completes the proof of (iii).
	\endproof

The next lemma, establishing the upper semicontinuity of the sunlight functional $\S$
w.r.t.~weak convergence of measures, provides the main ingredient 
in the proof of existence of optimal measures. 
We recall that the weak convergence of measures $\mu_k\wto \mu$ means
\bel{wlim}\lim_{k\to\infty} \int \vp \, d\mu_k~=~\int \vp \, d\mu\qquad\qquad
\hbox{for every } ~~\vp\in \C^0_c(\R^d).\eeq

In the following we consider a sequence of positive Radon measures $(\mu_k)_{k\geq 1}$, 
on $\R^d$, satisfying the usual assumptions

{\bf 1 - Boundedness:}  there exists a constant $C$ such that 
\bel{mkb}
\mu_k(\R^d)\leq C\qquad\hbox{for all}~~k\geq 1.\eeq

{\bf  2 - Tightness:} for every $\ve>0$ there exists a radius $R_\ve$ such that
\bel{tight}\mu_k\Big(\bigl\{ x\in \R^d\,;~|x|>R_\ve\bigr\}\Big)~<~\ve.\eeq
\v
By a well known compactness theorem \cite{AGS, Bi}, this implies the existence of a weakly convergent 
subsequence: $\mu_{k_j} \wto \mu$.   
\v
\begin{lemma}\label{lemma4} 
Consider a weakly convergent
sequence of measures $\mu_k\wto \mu$, satisfying the boundedness and tightness conditions (\ref{mkb})-(\ref{tight}).    Then, for any unit vector $\bfn$ and 
every positive measure $\nu$, absolutely
continuous w.r.t.~Lebesgue measure on $\R^d$, one has
\bel{Susc}\S^\bfn(\mu;\nu)~\geq~\limsup_{k\to\infty} \S^\bfn(\mu_k;\nu).\eeq
\end{lemma}
\v
{\bf Proof.}  {\bf 1.} We start with the basic case where $\nu=0$ and all measures
$\mu_k$ are supported inside a ball $B(0, R)\subset\R^d$.

From the assumption it follows the weak convergence
$\mu^\bfn_k\wto\mu^\bfn$ of the projected measures.
Call 
$\Phi^\bfn_k$, $\Phi^\bfn$ respectively the density of the absolutely continuous
part of $\mu^\bfn_k$ and $\mu^\bfn$ w.r.t.~$(d-1)$-dimensional Lebesgue measure
on $E^\perp_\bfn$.   

Let $\ve>0$ be given.
According to the ``biting lemma" \cite{BM,  BC80, GMS},
 there exists 
a set $V_\ve\subset B(0,R)\subset E^\perp_\bfn$, with
\bel{mease}meas(V_\ve)~<~\ve,\eeq
and such that the following holds.
Let $\hat \mu_k^\bfn$ be the absolutely continuous measure 
on $E^\perp_\bfn$ whose density (w.r.t.~Lebesgue measure) is 
$$\Hat \Phi^\bfn_k(y)~=~\left\{ \bega{cl} \Phi^\bfn_k(y)\quad &\hbox{if}~~y\in B(0,R)\setminus V_\ve\,,\\[3mm]
0\quad &\hbox{otherwise.}
\enda\right.$$
Then, by possibly extracting a subsequence, we have the weak convergence
$$\hat \mu_k^\bfn~\wto~\hat \mu^\bfn,\qquad\qquad \Hat \Phi^\bfn_k~\wto~\Hat \Phi^\bfn\,,$$
Here the second arrow denotes weak convergence in $\L^1$.
Moreover, $\hat \mu^\bfn$ is the absolutely continuous measure having density 
$\Hat \Phi^\bfn$ w.r.t.~$(d-1)$-dimensional Lebesgue measure.
By (\ref{mease}) one has the obvious estimate
\bel{oe}
\int_{V_\ve} \Big(1- \exp \bigl\{ 
-\Phi^\bfn(y)\bigr\}\Big)\, dy~\leq~meas(V_\ve)~<~\ve.\eeq
Since $\mu^\bfn\geq \hat\mu^\bfn$, by (\ref{SSn}) and (\ref{oe}) 
 the total sunshine captured by the measure $\mu$ can now be estimated as
\bel{WCS}\bega{rl}
\S^\bfn(\mu)&\ds\geq~\int_{E^\perp_\bfn} \Big(1- \exp \bigl\{ 
-\Hat\Phi^\bfn(y)\bigr\}\Big)\, dy\\[4mm]
&\geq~\ds \limsup_{k\to \infty} \int_{E^\perp_\bfn} \Big(1- \exp \bigl\{ 
-\Hat\Phi_k^\bfn(y)\bigr\}\Big)\, dy
\\[4mm]
&=~\ds \limsup_{k\to \infty} \left(\int_{E^\perp_\bfn} \Big(1- \exp \bigl\{ 
-\Phi_k^\bfn(y)\bigr\}\Big)\, dy -\int_{V_\ve} \Big(1- \exp \bigl\{ 
-\Phi_k^\bfn(y)\bigr\}\Big)\, dy\right)\\[4mm]
&\geq~\ds\limsup_{k\to \infty}\, \S^\bfn(\mu_k)-\ve.
\enda
\eeq
Notice that  the concavity of the 
function $x\mapsto (1-e^{-x})$ was here used in the estimate of the weak limit.
Since $\ve>0$ was arbitrary, this proves the lemma in the basic case.
\v
{\bf 2.} Next, we still assume that the measures $\mu_k$ have uniformly bounded
support, say 
\bel{suppk}Supp(\mu_k)~\subseteq~ B(0,R)\qquad\qquad\forall k\geq 1\,,\eeq
but we allow the presence of an additional  positive measure $\nu$, having density
$g\in \L^1_{loc}(\R^d)$ w.r.t.~Lebesgue measure.
In the following we consider the cylinder
\bel{Gdef} 
\Gamma_{R}~\doteq~ \bigl\{y+t\bfn\,;~~y\in E^\perp_\bfn\,,~~r\in\R\,,~~ |y|\leq R,~|t| \leq R 
\bigr\}.\eeq
Let $\ve_0>0$ be given.  Then there exists $\rho_0>0$ such that
\bel{ig0}
\int_V g(x)\, dx~\leq~\ve_0\eeq
for every set
$V\subseteq \Gamma_R$ such that $ meas(V)\leq\rho_0$.
Calling $\omega_{d-1}$ the volume of the unit ball in $\R^{d-1}$, we choose 
$\rho_1>0$ so that 
\bel{ir1} \omega_{d-1} R^{d-1} \rho_1~<~\rho_0\,.\eeq
Then we choose 
$$-\infty~=~t_0~<~t_1~<~t_2~<~\cdots~<~t_N~<~t_{N+1}~=~+\infty$$
such that 
\bel{tp1}t_1<-R,\qquad t_N> R,\qquad  t_j-t_{j-1}< \rho_1\qquad\forall j=2,\ldots,N,\eeq
\bel{tp2}\mu\Big(\{ x\in\R^d\,;~~\la x,\bfn\ra = t_j\}\Big)~=~0\qquad\forall j=1,2,\ldots,N.\eeq
\v
{\bf 3.} Call $\mu^j_k~=~\chi_{\strut\{ \langle x,\bfn\rangle \geq t_j\}}\cdot\mu_k$
the restriction of the measure $\mu_k$ to the set where $\langle x,\bfn\rangle \geq t_j $, 
and let
$\mu^{\bfn,j}_k$, $\mu^{\bfn,j}$  the projections of $\mu^{j}_k$, $\mu^{j}$ on $E^\perp_\bfn$,
as in (\ref{mupro}).  Moreover, call $\Phi^{\bfn,j}_k$, $\Phi^{\bfn, j}$ the densities
of the absolutely continuous parts of $\mu^{\bfn,j}_k$, $\mu^{\bfn,j}$ w.r.t.~the $(d-1)$-dimensional
Lebesgue measure on $E^\perp_\bfn$.
 
The weak convergence $\mu_k\wto\mu$, together with the assumption (\ref{tp2}) implies the weak convergence
\bel{wjk}  \mu^{\bfn,j}_k~\rightharpoonup~\mu^{\bfn,j} \qquad\qquad\forall ~j=0,1,\ldots,N.
\eeq
Using again the  ``biting lemma" \cite{BC80}, 
we can find a set $V\subseteq B(0,R)\subseteq E_{\bfn}^\perp$, with 
\bel{meas0}\meas(V)~<~\ve_0,\qquad\qquad
(t_N-t_1)\cdot \meas(V)~<~\rho_0,
\eeq
and such that the following holds. 
Let $\Hat \mu_k^{\bfn,j}$ be the absolutely continuous measure on $E_\bfn^\perp$ whose density  is 
\bel{HPdf}\Hat \Phi^{\bfn,j}_k(y)~=~\left\{ \bega{cl} \Phi^{\bfn,j}_k(y)\quad &\hbox{if}~~y\in B(0,R)\setminus V,\\[3mm]
0\quad &\hbox{otherwise.}
\enda\right.\eeq
Then, by possibly extracting a subsequence, for every $j=1,\ldots,N$ we have the weak convergence
$$\Hat\mu^{\bfn,j}_k\rightharpoonup\Hat\mu^{\bfn,j},\qquad
\qquad\Hat\Phi^{\bfn,j}_k\rightharpoonup\Hat\Phi^{\bfn,j}\,.$$
Here the second arrow denotes weak convergence in $\mathbf{L}^1$. 
Moreover, $\Hat\mu^{\bfn,j}$ is the absolutely continuous measure 
on $E^\perp_\bfn$ with density $\Hat\Phi^{\bfn,j}$. 
\v
{\bf 4.}
For each fixed $y\in E^\perp_\bfn$,  the last integral in 
(\ref{Imny}) can be estimated from above and from below in terms of Riemann sums.
More precisely, for a given measure $\mu$, call
$$\mu^{j}~=~\chi_{\strut\{ \langle x,\bfn\rangle \geq t_j\}}\cdot\mu$$
the restriction of $\mu$ to the set $ \{ x\in\R^d\,;~~\langle x,\bfn\rangle \geq t_j\}$.
Let $\mu^{\bfn,j}$ be the projection of $\mu^j$ on $E^\perp_\bfn$, and let
 $\Phi^{\bfn,j}$ be the density of the absolutely continuous part of $\mu^{\bfn,j}$.
Since $\mu,\nu$ are both positive measures, one has
\bel{bes}\bega{l}\ds L(y)~\doteq~\sum_{j=2}^N
\int_{t_{j-1}}^{t_j} \left( g(y+s\bfn) \exp\left\{
-\int_s^{+\infty} g(y+t\bfn)\, dt\right\} 
\cdot \exp\Big\{ -\Phi^{\bfn, j-1}(y)\Big\}
\right)\, ds
\\[4mm]
\ds\leq~\int_{t_1}^{t_N} \left( g(y+s\bfn) \exp\left\{
-\int_s^{+\infty} g(y+t\bfn)\, dt\right\} 
\cdot \exp\Big\{ -\Phi^\bfn(y)\cdot\mu^y\bigl( [s, +\infty[ \,
\bigr)\Big\}
\right)\, ds\\[4mm]
\ds\leq~\sum_{j=2}^N
\int_{t_{j-1}}^{t_j}  \left( g(y+s\bfn) \exp\left\{
-\int_s^{+\infty} g(y+t\bfn)\, dt\right\} 
\cdot \exp\Big\{ -\Phi^{\bfn,j}(y)\Big\}
\right)\, ds~\doteq~U(y)\,.
\enda
\eeq
The difference between the upper and lower Riemann sums, 
on the right and the left hand side of (\ref{bes}), can be estimated by
\bel{U-L}\bega{l}\ds 
U(y)-L(y)~\leq~\sum_{j=2}^N
\int_{t_{j-1}}^{t_j} g(y+s\bfn)\, ds\cdot \Big[ \exp\bigl\{ -\Phi^{\bfn,j}(y)\bigr\}-\exp\bigl\{ -\Phi^{\bfn,j-1}(y)\bigr\}
\Big]\\[4mm]\ds
\qquad \leq~\max_j\left(
\int_{t_{j-1}}^{t_j} g(y+s\bfn)\, ds\right)\cdot \sum_{j=2}^N
\Big[ \exp\bigl\{ -\Phi^{\bfn,j}(y)\bigr\}-\exp\bigl\{ -\Phi^{\bfn,j-1}(y)\bigr\}
\Big] \\[4mm]\qquad\ds \leq~\max_j\left(
\int_{t_{j-1}}^{t_j} g(y+s\bfn)\, ds\right).\enda
\eeq
By (\ref{U-L}) and the choice of the points $t_j$ it now follows
\bel{ULe}\int_{|y|<R} \Big(U(y)-L(y)\Big)\, dy~\leq~
 \int_{|y|<R}\left( \sup_{2\leq j\leq N}  \int_{t_{j-1}}^{t_j} g(y+s\bfn)\, ds\right)   dy~\leq~\ve_0\,.\eeq
 Indeed, to prove the last inequality, consider a measurable subset $\Gamma^\sharp\subset\Gamma$ such that
$$\Gamma^\sharp~=~\Big\{ (y+s\bfn)\,;~~~|y|<R,~~s\in [t_{j(y)-1}, t_{j(y)}]\Big\},$$
where, for a.e.~$y\in B(0,R)$,
$$ \int_{t_{j(y)-1}}^{t_{j(y)}} g(y+s\bfn)\, ds~=~\max_{2\leq j\leq N} 
 \int_{t_{j-1}}^{t_j} g(y+s\bfn)\, ds.$$
 By (\ref{ir1}), 
 $$meas(\Gamma^\sharp)\leq ~\omega_{d-1}\cdot R^{d-1}\, \rho_1~<~\rho_0\,.$$
 Hence, by (\ref{ig0}) and the definition of $\Gamma^\sharp$,
 \bel{UL3}
 \int_{|y|<R}\left( \sup_{2\leq j\leq N}  \int_{t_{j-1}}^{t_j} g(y+s\bfn)\, ds\right)   dy~\leq~
 \int_{\Gamma^\sharp} g(x)\, dx~\leq~\ve_0\,.\eeq
\v
{\bf 5.} Recalling (\ref{Imny})-(\ref{Imn}) and using  (\ref{U-L}), we obtain
\bel{ME1}\bega{l}\ds
\S^\bfn(\mu_k,\nu)-\S^\bfn(\mu,\nu)~=~\int_{E^\perp_\bfn} \Big( \S^\bfn_{\mu_k,\nu}(y) -\S^\bfn_{\mu,\nu}(y)
\Big)\, dy \\[4mm]
\leq~\ds
\int_{|y|<R} \exp\left\{ -\int g(y+s\bfn)\, dt\right\}\cdot\Big(\exp\{ -\Phi^\bfn(y)\} - \exp\{ -\Phi^\bfn_k(y)\}\Big)\, dy
\\[4mm]\ds
\quad +\int_{|y|<R}\sum_{j=1}^N\int_{t_{j-1}}^{t_j}
g(y+s\bfn) \exp\left\{ -\int_s^{+\infty}
 g(y+t\bfn)\, dt\right\}\, ds\\[4mm]
 \qquad\qquad\qquad\qquad \ds \cdot\Big(\exp\{ -\Phi^{\bfn,j-1}(y)\} - \exp\{ -\Phi_k^{\bfn,j}(y)\}\Big)\, dy
 \\[4mm]
 \leq~\ds
\int_{|y|<R} \exp\left\{ -\int g(y+s\bfn)\, dt\right\}\cdot\Big(\exp\{ -\Phi^\bfn(y)\} - \exp\{ -\Phi^\bfn_k(y)\}\Big)\, dy
\\[4mm]\ds
\quad +\int_{|y|<R}\sum_{j=1}^N\int_{t_{j-1}}^{t_j}
g(y+s\bfn)  \exp\left\{ -\int_s^{+\infty}
 g(y+t\bfn)\, dt\right\}\, ds\\[4mm]
 \qquad\qquad\qquad\qquad \ds \cdot\Big(\exp\{ -\Phi^{\bfn,j}(y)\} - \exp\{ -\Phi_k^{\bfn,j}(y)\}\Big)\, dy\\[4mm]\ds
\quad 
+\int_{|y|<R}\bigl( U(y)-L(y)\bigr)\, dy\\[4mm]
=~I_{1,k}+I_{2,k}+I_3\,.
\enda
\eeq
As $k\to\infty$,
the limits of the first two  integrals can be estimated as in Step {\bf 1}.  Indeed, recalling 
the properties (\ref{meas0})-(\ref{HPdf}) we obtain
\bel{I2}\limsup_{k\to\infty}\int_{|y|<R,~ y\notin V} 
\exp\left\{ -\int g(y+s\bfn)\, dt\right\}\cdot\Big(\exp\{ -\Hat \Phi^\bfn(y)\} - \exp\{ -\Hat\Phi^\bfn_k(y)\}\Big)\, dy
~\leq~0,\eeq
\bel{I3j}
\bega{l}\ds \limsup_{k\to\infty}\int_{|y|<R,~ y\notin V} \int_{t_{j-1}}^{t_j}
g(y+s\bfn)  \exp\left\{ -\int_s^{+\infty}
 g(y+t\bfn)\, dt\right\}\, ds \\[4mm]
 \ds\qquad\qquad \qquad\qquad \qquad\qquad \cdot\Big(\exp\{ -\Hat\Phi^{\bfn,j}(y)\} - 
 \exp\{ -\Hat\Phi_k^{\bfn,j}(y)\}\Big)\, dy\
~\leq~0.\enda \eeq
Moreover, by (\ref{ULe}) we already know that 
$I_3\leq\ve_0$.
{}From (\ref{meas0}) and the above inequalities we conclude 
\bel{sc1}
\limsup_{k\to\infty} \,(I_{1,k}+I_{2,k})+I_3~\leq~\ve_0 + \ve_0\,.\eeq
Since $\ve_0>0$ was arbitrary, this proves (\ref{Susc}) in the case where
the supports of the measures $\mu_k$ are uniformly bounded.
 \v
{\bf 6.} Finally, using the tightness assumption (\ref{tight}), we remove the assumption that the measures $\mu_k$ have uniformly bounded 
support.

For any given $\ve>0$, by (\ref{tight}) there exists  a radius $R$ sufficient large such that
$$\mu_k\Big(\{x\in \R^d;~|x|> R\}\Big)~<~\ve$$
for every $k\geq 1$. 
Without loss of generality, we can assume that 
\bel{wko}\mu\Big(\{x\in \R^d;~|x|= R\}\Big)~=~0.\eeq
Calling $B_R$ the open ball centered at the origin with radius $R$, we denote by
 $\mu_k^\flat, \mu_k^\sharp$ the restrictions of $\mu_k$ to  $B_R$ and $\R^d\setminus B_R$, respectively. 
 The measures $\mu^\flat$, $\mu^\sharp$ are defined similarly.
 By the weak convergence $\mu_k\wto\mu$ together with 
 (\ref{wko}) it follows the weak convergence $\mu_k^\flat\rightharpoonup \mu^\flat$. By Lemma~\ref{l:sl2}, for every $k$
 one has
$$\S^\bfn(\mu_k^\flat,\nu)~\geq~\S^\bfn(\mu_k,\nu)-\mu_k^\sharp\Big(\R^d\setminus B_R\Big)~\geq~\S^\bfn(\mu_k,\nu)-\ve\,.$$
Since the measures $\mu^\flat_k$ have uniformly bounded support, 
 by the previous analysis we conclude
$$\S^\bfn(\mu,\nu)~\geq~\S^\bfn(\mu^\flat,\nu)~\geq~\limsup_{k\to\infty}\S^\bfn(\mu_k^\flat,\nu)~\geq~\limsup_{k\to\infty} \S^\bfn(\mu_k,\nu)-\ve.$$
Since $\ve$ is arbitrary, this completes the proof.
\endproof

{}From the above lemma one easily obtains the upper semicontinuity of
the functional in (\ref{SSe}).

\begin{lemma}\label{lemma5}
Consider a weakly convergent
sequence of measures $\mu_k\wto \mu$, satisfying the boundedness and tightness conditions (\ref{mkb})-(\ref{tight}).    Then, for any positive, integrable function $\eta\in \L^1(S^{d-1})$ and 
every positive measure $\nu$, absolutely
continuous w.r.t.~Lebesgue measure on $\R^d$, one has
\bel{usc}\S^\eta(\mu;\nu)~\geq~\limsup_{k\to\infty} \S^\eta(\mu_k;\nu).\eeq
\end{lemma}
\v
{\bf Proof.}   By the boundedness assumption (\ref{mkb}) and the estimate (i) in Lemma~\ref{l:sl2}, for each $k\geq 1$
we have 
\bel{es6}\S^\bfn(\mu_k;\nu)~\leq~\mu_k(\R^d)~\leq~ C.\eeq
This implies
\bel{es7}\bega{l}\ds
\limsup_{k\to\infty} \S^\eta(\mu_k,\nu)~=~\limsup_{k\to\infty}
\int_{S^{d-1}}\eta(\bfn) \S^\bfn(\mu_k,\nu)~d\bfn~\leq~\int_{S^{d-1}}
\limsup_{k\to\infty} \eta(\bfn) \S^\bfn(\mu_k,\nu)~d\bfn\\[4mm]
\ds\qquad\leq~\int_{S^{d-1}}\eta(\bfn) \S^\bfn(\mu,\nu)~d\bfn
~=~\S^\eta(\mu,\nu).\enda\eeq
Here the first inequality is valid because, by (\ref{es6}),  all  integrand functions are pointwise bounded by
the function
$C\,\eta(\cdot)\in  \L^1(S^{d-1})$.  The second inequality follows from Lemma~\ref{lemma4}.
\endproof
\v

\section{Harvest functionals}
\label{s:harvest}
\setcounter{equation}{0}
We now consider a utility functional 
associated with roots, whose the main goal is to collect
moisture and nutrients from the ground. 
To model the efficiency of a root, consider a scalar function  $u(\cdot)$
and a positive measure
$\mu$.
We think of  $u(x)$ as the density 
of water+nutrients at the point $x$, while $\mu$ is the
density of root hair cells, which absorb fluids from the soil.
Since these fluids diffuse through the soil and are
harvested by the root, one expects that  $u$ satisfies a parabolic
equation of the form
\bel{pare}u_t ~=~\Delta u + f(x,u)-u\,\mu.\eeq
Since we are interested in average values over long periods of time, 
we look at the equilibrium states for (\ref{pare}).   
Throughout the following, we assume that
\begi
\item[(A1)]  ~$\Omega\subset\R^d$ is a bounded, connected open set with $\C^2$ boundary.
\item[(A2)]~$f:\ov\Omega\times\R\mapsto\R$ is a bounded, continuous function such that, for 
some constants $M,K$,
\bel{fprop}
f(x,0)~\geq~0,\qquad f(x,M)~\leq~0,\qquad |f(x,u)|~\leq~K,\qquad\qquad\forall x\in \ov\Omega,~~u\in [0,M]\,.\eeq
\item[(A3)] $\mu$ is a positive Radon measure supported on the compact set $\ov\Omega$.
\endi
We consider solutions $u:\ov\Omega\mapsto [0,M]$ of the elliptic problem with measure-valued coefficients
\bel{emu}
\Delta u + f(x,u) - u\,\mu~=~0,\eeq
and  Neumann boundary conditions
\bel{NBC} \partial_{\bfn(x)} u(x)~=~0\qquad\qquad x\in\partial\Omega.\eeq
Here $\bfn(x)$  denotes the unit outer normal vector
at the boundary point $x\in\partial\Omega$, while
$\partial_\bfn u= \bfn(x)\cdot \nabla u(x)$ is the derivative of $u$ in the
normal direction.

In alternative,  we shall also  consider  Dirichlet
 boundary conditions
\bel{DBC} u(x)~=~0\qquad\qquad x\in \partial\Omega.\eeq

Observe that, if the measure $\mu$ has a smooth density $h(\cdot)$ w.r.t.~Lebesgue measure,
then the equation (\ref{emu}) takes the form
\bel{ehu}
\Delta u + f(x,u) - h(x) \,u~=~0.\eeq
By the assumption (A2), the constant function $u_*(x)=0$ is a subsolution, while
$u^*(x)=M$ is a supersolution.   A standard comparison argument  now implies that 
the semilinear elliptic problem
(\ref{ehu}), (\ref{NBC}) has at least one solution $u:\ov\Omega\mapsto [0,M]$.
\v
Elliptic problems with measure data have been 
studied in several papers \cite{BG1, BGO, DMOP}
and are now fairly well understood. 
A key fact is that, roughly speaking, the Laplace operator ``does not see" sets with zero capacity.
Following \cite{BG1, BGO} we thus call $\M_b$ the set of all bounded Radon measures on $\ov\Omega$.
Moreover, we denote by $\M_0\subset\M_b$ the family of measures which vanish on Borel sets with zero capacity,
so that
\bel{mu0}\hbox{cap}_2(V)~=~0\qquad\implies\qquad \mu(V)~=~0.\eeq
For the definition and basic properties of capacity we refer to \cite{EG}.
Every measure $\mu\in\M_b$ can be uniquely decomposed as a sum
\bel{mudec}\mu~=~\mu_0 + \mu_s,\eeq
where $\mu_0\in \M_0$ while the measure $\mu_s$ is supported on a set with zero capacity.
In the definition of solutions to (\ref{emu}), the presence of the singular measure $\mu_s$ is disregarded.

\begin{remark}  {\rm  If $\mu$ is an arbitrary Radon measure
and $u$ is a measurable  function defined up to a set of zero Lebesgue measure, 
 the product $u\mu$ may not be well defined.
In the present setting, however, we claim that the product measure $u\,\mu_0$ 
is uniquely defined.  
Indeed,
calling
 $$\avint_V u\, dx~=~{1\over \hbox{meas}(V)}\,\int_V u\, dx$$
the average value of $u$ on a set $V$, for each $x\in \ov\Omega$ we 
can consider the limit
\bel{upv}
u(x)~=~\lim_{r\downarrow 0}~\avint_{\Omega\cap B(x,r)} u(y)\, dy.\eeq
As proved in  \cite{FZ}, if $u\in H^1(\Omega)$ then
the above limit exists at all points $x\in \ov\Omega$ with the possible exception of a set 
whose capacity is zero.   Since $\mu\in \M_0$, we conclude that the measure $u\mu_0$ is well defined.
}\end{remark}
 \v
\begin{definition} Let $\mu$ be a measure in $\M_b$,  decomposed as in (\ref{mudec}).   
 \begi
\item[(i)] A function
$ u\in\L^\infty(\Omega)\cap H^1(\Omega)$, with pointwise values given by (\ref{upv}), is a solution to  
the elliptic problem (\ref{emu})-(\ref{NBC}) if
 \bel{nsol}
 -\int_\Omega \nabla u\cdot \nabla\vp\, dx+\int_\Omega f(x,u)\vp\, dx -\int_{\ov\Omega} u\vp\, d\mu_0~=~0\eeq
 for every test function $\vp\in \C^\infty_c(\R^d)$.
\item[(ii)] 
A function
$ u\in\L^\infty(\Omega)\cap H^1_0(\Omega)$, with pointwise values given by (\ref{upv}), is a solution to  
the elliptic problem (\ref{emu}), (\ref{DBC}) if
 \bel{dsol}
 -\int_\Omega \nabla u\cdot \nabla\vp\, dx+\int_\Omega f(x,u)\vp\, dx -\int_\Omega u\vp\, d\mu_0~=~0\eeq
 for every test function $\vp\in \C^\infty_c(\Omega)$.
 \endi
\end{definition}
\v
We can now state the
 main existence result for solutions to (\ref{emu}). The proof closely follows the arguments in \cite{BCS}.
\v
\begin{thm}\label{thm3.1} Under the assumptions
(A1)--(A3), the elliptic problem   (\ref{emu}) with Neumann boundary conditions (\ref{NBC}) has at least one solution $u:\Omega\mapsto[0,M]$.
The same is true in the case of  Dirichlet boundary conditions (\ref{DBC}).
 \end{thm}
 \v
 {\bf Proof.}  
 Without loss of generality, we can assume that $\mu=\mu_0$, so that (\ref{mu0}) holds.
 \v
 {\bf 1.} We first consider  the case of Neumann boundary conditions. 
 Let $\Omega_\ve=\{ x\in\R^d\,;~~d(x,\ov\Omega)<\ve\}$ be a neighborhood
of radius $\ve>0$ around the compact set $\ov\Omega$. Following \cite{E}, 
we can construct a bounded, linear extension operator 
$E:H^1(\Omega)\mapsto H^1_0(\Omega_\ve)$.

Since $\mu$ is a bounded Radon measure on $\Omega_\ve$ which vanishes on sets of zero capacity,
by the analysis in \cite{DMOP} it follows that $\mu\in \L^1(\Omega_\ve)\oplus H^{-1}(\Omega_\ve)$. 
More precisely, there exist functions 
$\phi_0\in\L^1(\Omega_\ve)$ and $\phi_1,\ldots\phi_d\in \L^2(\Omega_\ve)$ such that 
\bel{murep}\int_{\Omega_\ve}  \vp\,d\mu~=~\int_{\Omega_\ve} \phi_0\, \vp\, dx-\sum_{i=1}^d \int_{\Omega_\ve}
\phi_i \,\vp_{x_i}\, dx\eeq
for every test function $\vp\in \C^\infty_c(\Omega_\ve)$.  Hence the same holds for every 
$\vp\in H^1_0(\Omega_\ve)$.
\v
{\bf 2.} By slightly shifting the measure $\mu$ in the interior of  the domain $\Omega$  and performing a mollification, we 
construct sequences of smooth functions
$\phi_{0,n}, \phi_{1,n},\ldots, \phi_{d,n}$ such that
\bel{phic}
\lim_{n\to\infty} \|\phi_{0,n}- \phi_0\|_{\L^1(\Omega_\ve)}~=~0,\qquad \lim_{n\to\infty}\|\phi_{j,n}- \phi_j\|_{\L^2(\Omega_\ve)}~=~0,\qquad\qquad j=1,\ldots,d.\eeq
 Moreover, the measures $\mu_n$ with density $h_n= \phi_{0,n} + \sum_j (\phi_{j,n})_{x_j}$ 
 w.r.t.~Lebesgue measure are nonnegative and supported in the interior of $\Omega$.
\v
{\bf 3.} By a standard comparison argument, for  each $n\geq 1$ we obtain the existence of a classical solution 
$u_n:\ov\Omega\mapsto [0,M]$
to the elliptic equation
\bel{een}
\Delta u+f(x,u) - h_n(x) u~=~0\qquad\qquad x\in \Omega\,,\eeq
with Neumann boundary conditions (\ref{NBC}).
Multiplying by $u_n$ and integrating by parts, one obtains
\bel{ees}
\int_\Omega
\Big[  \Delta u_n (x) +f(x, u_n(x))- h_n(x)u_n(x)\Big]  u_n(x)\, dx~=~0\eeq
Recalling that $h_n\geq 0$, $|f|\leq K$, and $u_n\in [0,M]$, we obtain 
\bel{h1es} \int_\Omega |\nabla u_n(x)|^2 \, dx~\leq~ \int_\Omega\bigl| f(x, u_n(x)) u_n (x)\bigr|\, dx ~\leq~
\hbox{meas}(\Omega)\cdot K\,M\,.\eeq
As a consequence, the norms $\|u_n\|_{H^1(\Omega)}$ remain uniformly bounded.
Therefore, the norms of the extensions $\|Eu_n\|_{H^1(\Omega_\ve)}$  are bounded as well.
\v
{\bf 4.} 
Thanks to the previous estimates, by possibly 
taking a subsequence and relabeling, we can assume the strong convergence
\bel{c4}\|u_n-u\|_{\L^2(\Omega)} ~\to ~0\eeq
and the weak convergence
\bel{c5}u_n\rightharpoonup u\qquad\qquad\hbox{in}~H^1(\Omega)\,,\eeq
\bel{c6}Eu_n\rightharpoonup Eu\qquad\qquad\hbox{in}~H^1_0(\Omega_\ve)\,,\eeq
for some function $u\in H^1(\Omega)$.
For every test function $\vp\in \C^\infty_c(\R^d)$ we now have
 \bel{uns}\bega{rl} 0&=~\ds
 \int_\Omega \Delta  u_n\, \vp\, dx+\int_\Omega f(x,u_n)\vp\, dx -\int_\Omega u_n\vp\, d\mu_n\\[4mm]
 &=~\ds
 -\int_\Omega \nabla u_n\cdot \nabla\vp\, dx+\int_\Omega f(x,u_n)\vp\, dx -\int_{\Omega_\ve} \left(\phi_{0,n} +\sum_{j=1}^d 
 (\phi_{j, n})_{x_j}\right)(Eu_n)\,\vp\, dx\\[4mm]
&=~\ds
 -\int_\Omega \nabla u_n\cdot \nabla\vp\, dx+\int_\Omega f(x,u_n)\vp\, dx -\int_{\Omega_\ve} \phi_{0,n} (Eu_n)\,\vp\, dx
 +\sum_{j=1}^d\int_{\Omega_\ve} 
 \phi_{j, n}(Eu_n)_{x_j}\,\vp\, dx\,.
 \enda
 \eeq
 Letting $n\to\infty$, by the strong convergence in (\ref{phic}) and (\ref{c4}) and  the weak convergence in
 (\ref{c5})-(\ref{c6}), we obtain
\bel{l1}\lim_{n\to\infty}\left( -\int_\Omega \nabla u_n\cdot \nabla\vp\, dx+\int_\Omega f(x,u_n)\vp\, dx \right)~=~
-\int_\Omega \nabla u\cdot \nabla\vp\, dx+\int_\Omega f(x,u)\vp\, dx\,,\eeq
 \bel{l2}\bega{l}\ds\lim_{n\to\infty}\left(\int_{\Omega_\ve} \phi_{0,n} (Eu_n)\,\vp\, dx
 +\sum_{j=1}^d\int_{\Omega_\ve} 
 \phi_{j, n}(Eu_n)_{x_j}\,\vp\, dx\right)\\[4mm]
 \ds\qquad\qquad ~=~\int_{\Omega_\ve} \phi_{0} (Eu)\,\vp\, dx
 +\sum_{j=1}^d\int_{\Omega_\ve} 
 \phi_{j}(Eu)_{x_j}\,\vp\, dx~=~\int_{\ov\Omega} u\vp\, d\mu\,.\enda\eeq
Together, (\ref{l1}) and (\ref{l2}) yield (\ref{nsol}), completing the proof in the case of Neumann boundary conditions.
\v
{\bf 5.} In the case of Dirichlet boundary conditions, without loss of generality we can assume that the measure $\mu$
is entirely supported in the interior of $\Omega$.  Indeed, since $u=0$ on the boundary, the part of $\mu$
supported on the boundary $\partial\Omega$ does not give any contribution to the right hand side of (\ref{dsol}).
We can thus use  the representation theorem in \cite{DMOP} directly on the set $\Omega$, and find
functions $\phi_0\in\L^1(\Omega)$, $\phi_1,\ldots,\phi_d\in \L^2(\Omega)$ such that
\bel{murep2}\int_\Omega  \vp\,d\mu~=~\int_\Omega \phi_0\, \vp\, dx-\sum_{i=1}^d \int_\Omega
\phi_i \,\vp_{x_i}\, dx\eeq
for every test function $\vp\in \C^\infty_c(\Omega)$. The proof is then 
achieved by  the same arguments as before.
\endproof 
 \v
{}From the proof of the above theorem, one can also obtain a comparison result.  As usual, we say that two Radon 
 measures satisfy $\tilde \mu\leq \mu$ if
 $\tilde\mu(V)\leq \mu(V)$ for every Borel set $V$.
 \begin{lemma}\label{l:comp}
 Let the assumptions (A1)-(A2) hold and consider two measures $\mu\geq\tilde \mu$, both satisfying (A3).
 Let $u:\ov\Omega\mapsto [0,M]$ be a solution of (\ref{emu})-(\ref{NBC}).     Then, replacing $\mu$ with $\tilde\mu$,
one can find a corresponding solution
 $\tilde u:\ov\Omega\mapsto [0,M]$ such that $\tilde u(x)\geq u(x)$ for every  $x$.
 
 The same result holds in the case of Dirichlet boundary conditions (\ref{DBC}).\end{lemma}
 \v
{\bf Proof.} Consider  the case of Neumann boundary conditions.   Define
the positive measure $\mu^* = \mu-\tilde\mu$.  Performing the construction described
in step {\bf 2} of the proof of Theorem~\ref{thm3.1}, with the same shifts and the same mollifications 
applied to all three measures, we obtain sequences of smooth functions 
$\phi_{i,n}$, $\tilde \phi_{i,n}$, $\phi_{i,n}^*$, for $i=0,1,\ldots, d$ and $n\geq 1$.
Since $\mu=\tilde\mu+\mu^*$, this implies that the corresponding densities of the mollified measures
satisfy
$$h_n(x)~=~\tilde h_n(x)+h_n^*(x)~\geq~\tilde h_n(x)~\geq~0.$$
Since $u_n$ is a solution to  (\ref{een}), it is a subsolution to 
\bel{subs} \Delta u+f(x,u) - \tilde h_n u~=~0, \eeq
always with Neumann boundary conditions (\ref{NBC}).
By a standard comparison argument,  there exists a solution
$\tilde u_n:\ov\Omega\mapsto [0,M]$ to (\ref{subs}), (\ref{NBC}).
such that 
$$\tilde u_n(x)~\geq ~u_n(x)\qquad\qquad  \hbox{for all}~~ x\in \ov\Omega.$$
By  taking  limits as $n\to \infty$, the result is proved.

The case of Dirichlet boundary conditions can be handled by the same technique.
\endproof
\v
We can now introduce a harvest functional, defined for solutions of 
 (\ref{emu}) with Neumann or Dirichlet boundary conditions.
 \v
\begin{definition}
Given a positive Radon measure $\mu$ on $\ov\Omega$ and a solution $u$ of (\ref{emu})-(\ref{NBC}), or 
(\ref{emu}), (\ref{DBC}), 
the {\bf total harvest}  is defined as
\bel{HF1}\H(u,\mu)~\doteq~\int_{\ov\Omega}~ u\, d\mu_0\,,\eeq
where $\mu=\mu_0+\mu_s$ is the decomposition introduced at (\ref{mudec}). \end{definition}
\v
In the case of Neumann boundary conditions, following  \cite{BCS} a more precise construction can be performed.
Let $G=G(t,x;y)$ be the Green function for the heat equation 
\bel{eel}\left\{
\bega{rl}
w_t ~=~\Delta w,&\qquad \qquad t>0\,,~~x\in\Omega\,,
\\
 \partial_{\bf n} w~=~0,&\qquad\qquad t>0\,,~~ x\in \partial \Omega\,.
\enda\right.
\end{equation}
As it is well known \cite{F}, 
for each fixed $y\in \Omega$  the function $G(\cdot, \cdot;y)$
provides a solution to (\ref{eel}) such that
\bel{intone}
\int_\Omega G(t,x;y)\, dx ~=~1\,,\qquad\qquad \lim_{t\downarrow
 0} \int_\Omega G(t,x;y) \, \phi(x)\, dx
~=~\phi(y)\eeq
for every $\phi\in \C(\ov\Omega)$.
The solution
of (\ref{eel}) with a continuous initial data $w(0,x)=\phi(x)$ is thus given by
$$w(t,x)~=~\int_\Omega G(t,x;y)\,\phi(y)\, dy\qquad\qquad   t>0\,,~~x\in \Omega\,.$$

Let now $u\geq 0$ be any function such that 
\bel{eep}
\bega{rll}
\Delta u &\geq~-K&\qquad  \hbox{on}~\Omega\,,
\\
 \partial_{\bf n} u&=~0&\qquad \hbox{on}~ \partial \Omega\,.
\enda
\end{equation}
In particular, $u$ could be the
solution to the elliptic problem (\ref{ees}) constructed in Theorem~\ref{thm3.1}.
For any $t>0$, consider the averaged function
\bel{avp}
u^{(t)}(x)~\doteq~\int_\Omega G(t, x;y)\, u(y)\,dy\,.\eeq
Using the boundary conditions in (\ref{eel}) and (\ref{eep}) to integrate by parts,
by the first equations in (\ref{intone}) and (\ref{eep})
one obtains
\bel{der}\bega{rl}\ds{d\over dt} u^{(t)}(x)&=~\ds{d\over dt} \int_\Omega G(t, x;y)\, u(y)\,dy
~=~\int_\Omega G_t(t, x;y)\, u(y)\,dy\cr & \cr
&\ds=~\int_\Omega \Delta G(t, x;y)\, u(y)\,dx~=~\int_\Omega 
G(t, x;y)\, \Delta u(y)\,dy~\geq~-K\,.\enda
\eeq
As a consequence, for every $x\in \Omega$ the map $t\mapsto u^{(t)}(x)+Kt$ is nondecreasing. 
Since every function  $x\mapsto u^{(t)}(x)$ is uniformly continuous on $\Omega$, it admits a 
continuous extension to the 
closure $\ov\Omega$.  At each $x\in \ov\Omega$
we can thus uniquely define the value $u(x)$ by setting
\bel{limp}u(x)~=~
\lim_{t\downarrow 0} u^{(t)}(x)~=~\inf_{t>0} \Big(u^{(t)}(x)+Kt\Big).\eeq
The representation (\ref{limp}) shows that $u$ is the infimum of a decreasing sequence of continuous functions.
Hence $u$  is upper semicontinuous.

\v
We conclude this section by observing that,
in the case of  Neumann 
 boundary conditions,  the harvest functional can be equivalently written as
 \bel{HF2}\H(u,\mu)~\doteq~\int_{\Omega}~ f(x,u(x))\, dx\,.\eeq
 In the case of Dirichlet
 boundary conditions, assuming that the solution is $\C^1$ in a neighborhood of 
 the boundary $\partial \Omega$, the harvest functional can be expressed as
  \bel{HF3}\H(u,\mu)~\doteq~\int_{\Omega}~ f(x,u(x))\, dx+\int_{\partial\Omega}\partial_{\bfn(x)} u(x)\, d\sigma\,.\eeq

\v
\section{Optimal irrigation patterns}
\label{s:7}
\setcounter{equation}{0}
This section provides a brief review of ramified transport and 
optimal irrigation.
To fix the ideas, throughout the following we assume
\begi
\item[(A4)]
 $\Omega\subset\R^d$ is a connected, open set with Lipschitz boundary, whose closure contains the origin:
${\bf 0}\in \ov \Omega$.  
  \endi  
  Given $\alpha\in [0,1]$,
to define the $\alpha$-irrigation cost  of a bounded, positive measure $\mu$ on $\ov\Omega$,
 we shall follow
the Lagrangian approach of 
Maddalena, Morel, and Solimini \cite{MMS}.

Let $\kappa=\mu(\ov\Omega)$ be the total mass to be transported and let 
$\Theta=[0,\kappa]$.
We think of each $\theta\in \Theta$ as a ``water particle".
A measurable map 
\bel{iplan}
\chi:\Theta\times [0,1]~\mapsto~ \ov\Omega\eeq
is called an {\bf admissible irrigation plan} for the measure $\mu$ on $\ov\Omega$
if 
\begi
\item[(i)] For a.e.~$\theta\in \Theta$, the map
$t\mapsto \chi(\theta,t)$ is Lipschitz continuous. 
\item[(ii)] At time $t=0$ all particles are at the origin:
$\chi(\theta,0)={\bf 0}\in\R^d$ for all $\theta\in\Theta$.
\item[(iii)] At time $t=1$ the push-forward of the Lebesgue measure on $[0,\kappa]$ through the map $\theta\mapsto 
\chi(\theta,1)$ coincides with the measure $\mu$.
In other words, for every open set $A\subset\R^d$ there holds
\bel{chi1}\mu(A)~=~\hbox{meas}\Big( \{ \theta\in \Theta\,;~~\chi(\theta,1)\in A\bigr\}\Big).\eeq
\endi
To the
 irrigation plan $\chi$ we now attach a cost $\E^\alpha$.  Toward this goal,
given a point $x\in \R^d$ we first compute 
how many paths go through the point $x$.  
This is described by
\bel{chi}|x|_\chi~=~\meas\Big(\bigl\{\theta\in \Theta\,;~~\chi(\theta,t)= x~~~\hbox{for some}~~t\in [0,1]\bigr\}\Big).\eeq
We think of $|x|_\chi$ as the {\it total flux going through the
point $x$}.
\v
\begin{definition}
For a given $0<\alpha\leq 1$,
the total cost of the irrigation plan $\chi$ is
\bel{TCg}
\E^\alpha(\chi)~\doteq~\int_\Theta\left(\int_0^1 \bigl|\chi(\theta,t)
\bigr|_\chi^{\alpha-1} \cdot |\chi_t(\theta,t)|\, dt\right)
d\theta\,.\eeq
If $\mu$ is a positive, bounded  Radon measure supported on $\ov \Omega$,
the  {\bf $\alpha$-irrigation cost} of $\mu$
is defined as 
\bel{Idef}\I^\alpha(\mu)~\doteq~\inf_\chi \E^\alpha(\chi),\eeq
where the infimum is taken over all admissible irrigation plans.
\end{definition}
\v
\begin{remark}  {\rm  
In the optimal irrigation problem, water has to be transported from a central well located at the origin 
${\bf 0}\in\R^d$ to various locations inside $\ov\Omega$.   We think of  $\chi(\theta,t)$ as the 
position of the water particle $\theta$ at time $t$. 
The factor $|\chi(\theta,t)|_\chi^{\alpha-1}$ models the assumption that water is transported through a network
of pipes,  whose cost is proportional to the product 
$$[\hbox{length}]\times[\hbox{flux}]^\alpha .$$
When 
$\alpha=1$ the integral in (\ref{TCg}) reduces to
$$\E^1(\chi)~\doteq~\int_\Theta\left(\int_0^1|\chi_t(\theta,t)|\, dt\right)
d\theta~=~\int_\Theta [\hbox{length of $\chi(\theta,\cdot)$}]\, d\theta
\,.$$
If $\Omega$ is convex, the minimum irrigation cost is trivially achieved by transporting each particle along a straight line,
hence
\bel{rpp23}\I^1(\mu)~=~\int |x|\, d\mu(x).\eeq
On the other hand, when
$\alpha<1$, it becomes convenient to lump together several paths into a unique large pipe, and the optimal
irrigation pattern can have a complicated structure.}
\end{remark}
\begin{remark}  {\rm 
In an irrigation pattern,
what matters are only the paths $\{\chi(\theta,t)\,;~~t\in [0,1]\}\subset\R^d$, not the time law 
with which these paths are traversed.
Indeed, for each $\theta$ we could take a 
smooth bijection $\tau^\theta:[0,1]\mapsto [0,1]$  and 
consider the time-reparameterized path 
$\tilde\chi(\theta,t) = \chi(\theta, \tau^\theta(t))$. 
Then the irrigation plan $\tilde\chi$ 
has exactly the same cost as $\chi$.}\end{remark}
\begin{remark}  {\rm 
As  suggested by intuition,  irrigation plans with minimum cost do not have loops.  Namely:
\bel{oplan1}\chi(\theta,t_1)~=~\chi(\theta, t_2)\qquad\implies\qquad \chi(\theta, t)~=~\chi(\theta,  t_1)\quad
\forall t\in [t_1,t_2].\eeq
A further, useful property of optimal irrigation plans is 
\bel{oplan2} \chi(\theta_1,t_1)~=~\chi(\theta_2,t_2)\qquad\implies\qquad \bigl\{\chi(\theta_1,t)\,;~~t\in [0, t_1]\bigr\}
~=~ \bigl\{\chi(\theta_2,t)\,;~~t\in [0, t_2]\bigr\}.\eeq
}\end{remark}
\v
For the basic theory of ramified transport we refer 
to \cite{BS, MMS, MS, X3, X15}, or to the monograph \cite{BCM}.
The next lemmas review the existence and some basic properties of the irrigation functional.
\v
\begin{lemma}\label{lexi}  Let $\Omega$ be a domain satisfying (A4), let $\alpha\in [0,1]$, 
and let $\mu$ be a bounded, 
positive measure on $\ov\Omega$.   If there exists an admissible irrigation pattern 
with finite cost $\E^\alpha(\chi)<+\infty$,
then  the measure $\mu$ admits an optimal irrigation plan.
\end{lemma}
For a proof, see Proposition 3.41 in  \cite{BCM}.  
\v
\begin{lemma}\label{lemma3}
Let $\Omega\subset\R^d$ satisfy the assumptions in (A4) and let $\mu$, 
$ \mu_1,\mu_2$ be bounded, positive measures on 
$\ov\Omega$. Then 
\bel{Ia1}
\I^\alpha(\mu)~\geq~\bigl[\mu(\ov\Omega)\bigr]^{\alpha-1}\,\int_{\ov\Omega} |x|\, d\mu\,,\eeq
 \bel{Ia2}\I^\alpha(\mu_1)~\leq~\I^\alpha(\mu_1+ \mu_2)~\leq~\I^\alpha(\mu_1)+\I^\alpha( \mu_2).\eeq
\end{lemma}
\v
{\bf Proof.}
{\bf 1.} The first inequality follows immediately from
$$\bega{rl}\ds
\E^\alpha(\chi)&\ds\doteq~\int_\Theta\left(\int_0^1 \bigl|\chi(\theta,t)
\bigr|_\chi^{\alpha-1} \cdot |\chi_t(\theta,t)|\, dt\right)
d\theta~\geq~\bigl[\mu(\ov\Omega)\bigr]^{\alpha-1}\int_\Theta\left(\int_0^1 |\chi_t(\theta,t)|\, dt\right)
d\theta\\[4mm]
\ds &\ds\geq~\bigl[\mu(\ov\Omega)\bigr]^{\alpha-1}\int_\Theta |\chi(\theta,1)|\, d\theta~=~\bigl[\mu(\ov\Omega)\bigr]^{\alpha-1}\int_{\ov\Omega} |x|\, d\mu(x)\,.\enda$$
\v
{\bf 2.} 
Next, for $i=1,2$ let
$\kappa_i~=~\mu_i(\ov \Omega)$ and let $\chi_i:[0,\kappa_i]\times [0,1]\mapsto\ov\Omega$ 
be an admissible irrigation plan
for $\mu_i$.
Then the map $\chi:[0, \kappa_1+\kappa_2]\times [0,1]\mapsto\ov\Omega$ defined by
$$\chi(\theta,t)~=~\left\{\bega{cl} \chi_1(\theta,t)\qquad &\hbox{if}\quad \theta\in [0,\kappa_1],\\[3mm]
\chi_2(\theta-\kappa_1,t)\qquad &\hbox{if}\quad \theta\in \,]\kappa_1, \, \kappa_1+\kappa_2],\enda
\right.
$$
is an admissible irrigation plan for $\mu_1+\mu_2$.  Its cost is
$$\bega{rl}\E^\alpha(\chi)&=\ds~\int_0^{\kappa_1+\kappa_2}
\left(\int_0^1|\chi(\theta,t)|_\chi^{\alpha-1}\cdot  |\chi_t(\theta,t)|\, dt\right)d\theta\\[4mm]
&\leq~\ds\int_0^{\kappa_1}
\left(\int_0^1|\chi(\theta,t)|_{\chi_{1}}^{\alpha-1}\cdot  |\chi_t(\theta,t)|\, dt\right)d\theta+
\int_{\kappa_1}^{\kappa_1+\kappa_2}
\left(\int_0^1|\chi(\theta,t)|_{\chi_{2}}^{\alpha-1}\cdot 
 |\chi_t(\theta,t)|\, dt\right)d\theta\\[4mm]
&=~\ds \E^\alpha(\chi_1)+\E^\alpha(\chi_2).
\enda
$$
This proves the second inequality in (\ref{Ia2}).
\v
To prove the first inequality we shall use the representation (see Proposition 4.8 in \cite{BCM})
\bel{H1rep}\E^\alpha(\chi)~=~
\int_\Theta\left(\int_0^1 \bigl|\chi(\theta,t)
\bigr|_\chi^{\alpha-1} \cdot |\chi_t(\theta,t)|\, dt\right)
d\theta~=~\int_{\R^d} |x|_\chi^\alpha\, d\H^1(x),\eeq
where $d\H^1$ denotes integration w.r.t.~the 1-dimensional Hausdorff measure.

Let $\chi:[0, \kappa_1+\kappa_2]\mapsto\ov\Omega$ be an admissible irrigation plan for 
$\mu_1+\mu_2$.  By possibly performing a measure-preserving transformation of the interval 
$\Theta=[0, \kappa_1+\kappa_2]$ into itself, we can assume that the map
$\chi_1:[0,\kappa_1]\times [0,1]\mapsto \ov\Omega$, obtained by restricting $\chi$ 
to the subdomain where $\theta\in [0,\kappa_1]$,
 is an admissible irrigation plan for $\mu_1$.
Using (\ref{H1rep}) we obtain the obvious estimate
$$\E^\alpha(\chi_1)~=~\int_{\R^d} |x|_{\chi_1}^\alpha\, d\H^1(x)~\leq~\int_{\R^d} |x|_{\chi}^\alpha\, d\H^1(x)
~=~\E^\alpha(\chi).$$
In other words, given any admissible irrigation plan for $\mu_1+\mu_2$, one can find an admissible irrigation plan 
for $\mu_1$ with smaller or equal cost.  This proves the first inequality in (\ref{Ia2}).
\endproof
\v

\section{Optimal shape of tree branches}
\label{s:4}
\setcounter{equation}{0}
Based on the functionals introduced in the previous sections,
 we now consider a constrained optimization problem
for a measure $\mu$ on $\R^d$, which we think as the distribution of 
 leaves on a tree.
The payoff will be the total amount of sunlight captured by the leaves. This 
will be supplemented by the cost of transporting nutrients from the base of the trunk, 
located at the origin ${\bf 0}\in\R^d$ to all leaves of the tree.

To  formulate this optimization problem, we consider:
\begi
\item[(i)] An open domain $ \Omega\subseteq\R^d$ with Lipschitz boundary, such that ${\bf 0}\in \ov\Omega$.
\item[(ii)] Constants 
$c, \kappa_0>0$, and 
an exponent $0<\alpha\leq 1$ such that
\bel{ad}
1-{1\over d-1} ~<~\alpha~\leq~1\,.\eeq
\item[(iii)] A non-negative, integrable  function $\eta:\S^{d-1}\mapsto
\R_+$, determining the intensity of light coming from various directions.
\item[(iv)] An absolutely continuous positive measure $\nu$, with  continuous density function $g:\R^d\mapsto\R_+$, describing the 
density of external vegetation. 
\endi
We then consider  the optimization problem
\bel{OP1}
\hbox{maximize:}\quad \S^\eta(\mu;\nu) - c
\I^\alpha(\mu).\eeq
subject to 
\bel{muR} Supp(\mu)\subseteq\ov\Omega,\qquad\qquad 
\mu(\ov\Omega)~\leq~\kappa_0\,.\eeq
Here $\S^\eta(\mu;\nu)$ is the {\it sunlight functional} introduced at
(\ref{Imny})--(\ref{Im}), while  
$\I^\alpha(\mu)$ is the minimum cost to $\alpha$-irrigate the measure $\mu$, 
defined at (\ref{Idef}).

\v
\begin{remark}  {\rm One can think of  (\ref{muR}) as  a constraint on 
the  size of the tree, i.e.~on the total amount of leaves.
Notice that  the inequality in (\ref{muR}) is essentially equivalent to
\bel{mu<1}
\mu(\ov\Omega)~=~\kappa_0\,.\eeq
Indeed, given a measure $\mu$ with total mass $<\kappa_0$,
we can always add to $\mu$ a Dirac mass at the origin, 
of size $\kappa_0-\mu(\R^d)$.  This would come at zero 
transportation cost, and zero additional payoff. } \end{remark}
\begin{remark}  {\rm If $\mu$ is supported on a set of dimension $< d-1$, then 
$\S^\bfn(\mu,\nu)=0$.  On the other hand, if (\ref{ad}) fails,
then $\I^\alpha(\mu)=+\infty$ for every measure $\mu$ whose support
is NOT contained in a set of dimension $\leq d-1$.
In this case, the above optimization problem would only have trivial solutions, where
the measure $\mu$ is a point mass at the origin.}
\end{remark}

Using the semicontinuity of the functionals $\S^\eta$ and $\I^\alpha$, and 
deriving suitable a priori estimates, we now prove
\v
\begin{thm}\label{thm1}
In the above setting (i)--(iv),  the 
constrained optimization problem (\ref{OP1})-(\ref{muR}) has at least one solution.
\end{thm}

{\bf Proof.} {\bf 1.}
By Lemma~\ref{l:sl} and the bound (\ref{muR}) it follows
\bel{ub}\S^\eta(\mu,\nu) -c\I^\alpha(\mu)~\leq~\S^\eta(\mu)~\leq~\int_{S^{d-1}}
 \eta(\bfn)\, d\bfn\cdot \sup\,\bigl\{ \S^\bfn(\mu)\,;~|\bfn|=1\bigr\}~\leq~ \|\eta\|_{\L^1}\,\kappa_0\,,\eeq
 showing that the  functional in (\ref{OP1}) has a finite upper bound. 
Hence there exists a  sequence of positive measures
$(\mu_n)_{n\geq 1}$, all satisfying the conditions in (\ref{muR}), and  such that
\bel{mumax}\lim_{n\to\infty}\Big( \S^\eta(\mu_n;\nu) - c
\I^\alpha(\mu_n)\Big)~=~\sup_\mu\Big\{ \S^\eta(\mu;\nu) - c
\I^\alpha(\mu)\Big\}.\eeq
The supremum on the right hand side is taken over all positive measures satisfying (\ref{muR}).
\v
{\bf 2.} We claim that it is not restrictive to assume that the measures
$\mu_n$ have uniformly bounded support. More precisely
\bel{Suppm}
Supp(\mu_n)~\subseteq~\ov B(0,r_0)\cap  \ov\Omega,\eeq
where
\bel{rdef}r_0~
=~{\kappa_0^{1-\alpha}\over c\,\alpha}\|\eta\|_{\L^1}\,.\eeq
Indeed,  each measure $\mu_n$ can written as a sum: $\mu_n=\Hat \mu_n+\mu^*_n$, where 
$\Hat \mu_n$ is supported inside the closed ball $\ov B(0,r_0)$,  while $\mu_n^*$ is supported outside this ball.

By Lemma~\ref{lexi},  for each $n\geq 1$ there exists an optimal irrigation  plan $\chi_n$, i.e., a minimizer
of the irrigation cost for the measure
 $\mu_n$.  By possibly performing a measure-preserving transformation of $\Theta_n=[0, \mu_n(\ov\Omega)]$ into itself,
 it is not restrictive to assume that 
 $$\chi_n(\theta,\cdot)~=~\left\{\bega{cl} \Hat \chi_n(\theta,\cdot)\qquad &\hbox{if}\quad \theta\in [0,\Hat \kappa_n],\\[3mm]
\chi_n^*(\theta-\Hat \kappa_n,\cdot)\qquad &\hbox{if}\quad \theta\in \,]\Hat \kappa_n, \, \kappa_n],\enda
\right.
$$
 where $\Hat \chi_n$ is an irrigation plan for $\Hat \mu_n$, while $\chi^*_n$ is an irrigation plan for $\mu^*_n$.
 By (\ref{muR}) we have
 $$\mu_n(\ov\Omega)~=~\Hat \mu_n (\ov\Omega)+\mu_n^*(\ov\Omega)
 ~\leq~
 \kappa_0\,.$$
  Since $\chi_n$ is optimal while $\Hat \chi_n$ is suboptimal,  
 the difference between the minimal  irrigation costs can be estimated as
\bel{dfic}
\begin{array}{rl}
\I(\mu_n)-\I(\Hat \mu_n)&\ds
\geq~\int_{\ov \Omega}\Big(|x|_{\Hat \chi_n}+|x|_{\chi^*_n}\Big)^\alpha\, d\H^1 - \int_{\ov\Omega}
 |x|_{\Hat \chi_n}^\alpha\,d\H^1\\[4mm]
&\ds
\geq~\int_{\ov \Omega}\alpha\,\kappa_0^{\alpha-1}\, |x|_{\chi^*_n}~d\H^1~
\geq~\alpha\,\kappa_0^{\alpha-1}\, r_0\cdot \mu^*_n(\ov\Omega).
\end{array}
\eeq
The second inequality comes from the fact that $(x+y)^\alpha-x^\alpha\geq \frac{\alpha}{\kappa_0^{1-\alpha}}y$, when $x\geq 0$, $y\geq 0$, and $ x+y\leq \kappa_0$.

On the other hand, by (\ref{qsy}) the difference in the sunlight functional can be estimated by
\bel{dfsl}
\S^\eta(\mu_n,\nu)-\S^\eta(\Hat \mu_n,\nu)~\leq~\|\eta\|_{\L^1}\cdot \mu^*_n(\ov\Omega).
\eeq
If the radius $r_0$ is chosen as in (\ref{rdef}), then
by (\ref{dfic})-(\ref{dfsl}) we have
\bel{Smf}\S^\eta(\Hat \mu_n,\nu) -c\I^\alpha(\Hat \mu_n)~\geq~\S^\eta(\mu_n,\nu) -c\I^\alpha(\mu_n)\eeq
for every $n\geq 1$.  By replacing each $\mu_n$ with $\Hat \mu_n$ we thus obtain a
maximizing sequence of measures whose supports are uniformly bounded.
\v
{\bf 3.} Thanks to the uniform boundedness of the supports,  by possibly taking a subsequence we can 
assume the weak convergence of measures: $\mu_n\rightharpoonup \bar\mu$, for some positive
measure $\bar\mu$ satisfying (\ref{muR}) as well.

By the lower semicontinuity of the irrigation cost $\I^\alpha$ (see Proposition 3.40 in \cite{BCM}),
it follows
$$\I^\alpha(\bar \mu)~\leq~\liminf_{n\to\infty} \I^\alpha(\mu_n).$$
On the other hand, the upper semicontinuity of the sunlight functional
$\mu\mapsto \S^\eta(\mu,\nu)$ proved in Lemma \ref{lemma5}
yields 
$$\S(\bar\mu,\nu)~\geq~\limsup_{n\to\infty} \S(\mu_n,\nu).$$
We conclude that $\bar \mu$ is an optimal solution to (\ref{OP1})-(\ref{mu<1}).
\endproof

\v
\section{Optimal shape of tree roots}
\label{s:5}
\setcounter{equation}{0}

In this section we consider constrained optimization problems
for a measure $\mu$ on $\R^d$, which we now think as the distribution of 
 root hair in the soil.  The payoff will be the total amount of water+nutrients collected
 by the roots. This will be supplemented by  the cost of transporting
 water   from the tips of the roots to the base of the  trunk.
 Under the same assumptions (A1)-(A2) in Section~\ref{s:harvest},
let  constants
$\alpha, c, \kappa_0>0$ be given, with 
\bel{ad2}
1-{1\over d-2}~<~\alpha~\leq~1\,.\eeq
We then consider  the optimization problem
\bel{OP2}
\hbox{maximize:}\quad \H(u,\mu) - c
\I^\alpha(\mu),\eeq
 among all positive measures $\mu$ on $\ov\Omega$ satisfying  the constraint 
\bel{muO}
\mu(\ov \Omega)~\leq~\kappa_0\,,\eeq
and all functions $u$ such that the couple 
$(u,\mu)$ provides a solution to the elliptic boundary value problem (\ref{emu}).
Here $\H(u,\mu)$ is the {\it harvest functional} introduced at
(\ref{HF1}), while  
$\I^\alpha(\mu)$ is the minimum cost to $\alpha$-irrigate the measure $\mu$, 
defined at (\ref{Idef}).

\v
\begin{remark}  {\rm 
 If $\mu$ is supported on a set of  zero capacity, then 
$\H(u,\mu)=0$.  As shown in chapter 5.9 of  \cite{AG},   
if  a set $A\subset\R^d$ has Hausdorff dimension $\leq d-2$, then its capacity is zero.

On the other hand (see \cite{BCM}), the minimum  irrigation cost $\I^\alpha(\mu)$  is bounded only if 
$\mu$ is supported on a set of dimension $< {1\over 1-\alpha}$.
To achieve a nontrivial solution of (\ref{OP2}), one thus needs the inequality   ${1\over 1-\alpha} > d-2$.
This motivates the condition (\ref{ad2}). }
\end{remark}
\v
Using the semicontinuity of the functionals $\H$ and $\I^\alpha$, we will prove the existence of optimal solutions.
We begin with the case of Neumann boundary conditions.
\v
\begin{thm}\label{thm2}   Let the assumptions (A1)-(A2) hold.
Then the 
maximization problem (\ref{OP2}), over all couples $(u,\mu)$ which satisfy (\ref{emu}), (\ref{NBC}), and (\ref{muO}), 
has an optimal solution.
\end{thm}
\v
{\bf Proof.} {\bf 1.} Call $\A$ the set of all admissible couples  $(u,\mu)$, satisfying (\ref{emu}), (\ref{NBC}), and (\ref{muO}).
Since every solution $u$ of (\ref{emu}) satisfies   $u(x)\in [0,M]$, 
 calling $\ov M$ the supremum over all admissible couples we have
\bel{ovM}\ov M~\doteq~\sup_{(u,\mu)\in\A} \bigl\{ \H(u,\mu)-\I^\alpha(\mu)\bigr\}~\leq~M\,\kappa_0\,.\eeq
Let $\{(u_n,\mu_n)\}_{n\ge 1}$ be a maximizing sequence.
 It is clearly not restrictive to assume that $\mu_n\in\M_0$
for every $n\geq 1$.

By (\ref{h1es})  we have the  bounds
$$\|u_n\|_{H^1}~\leq~C,\qquad\qquad  u_n(x)\in [0,M],$$
for some constant $C$ and every $n\geq 1$.
As remarked in Section~\ref{s:harvest}, the functions $u_n$ can be uniquely defined at every point $x\in \ov\Omega$
in terms of the limit (\ref{limp}).  Consider
the sequence of measures $\nu_n\doteq u_n\,\mu_n$.
By possibly taking a subsequence and relabeling we can assume
\bel{wklmt}\left\{\bega{l}
\ds~\nu_n\rightharpoonup \nu,\qquad\mu_n\rightharpoonup \mu\qquad\textrm{in the sense of weak convergence of measures,}\\[4mm]
\ds~~~~~~~~u_n\to u\qquad\qquad\qquad\textrm{strongly in $\L^2(\Omega)$ and a.e. in $\Omega$,}\\[4mm]
\ds~~~~~~~~u_n\rightharpoonup u\qquad\qquad\qquad\textrm{weakly in $H^1(\Omega)$.}
\enda\right.
\eeq
In addition, by Ascoli's theorem we can assume that, for every fixed $t>0$,
\bel{Gav}u_n^{(t)}(x)~=~\int_{\Omega} G(t,x,y)u_n(y)\, dy~\to~\int_{\Omega} G(t,x,y)u(y)\, dy~=~u^{(t)}(x)\eeq
as $n\to\infty$, uniformly for  $x\in\ov\Omega$.  Indeed, by choosing a subsequence we can 
achieve the convergence
in (\ref{Gav}) for every rational  $t>0$.  By continuity, this same subsequence satisfies
(\ref{Gav}) for every $t>0$.

\v
{\bf 2.} We claim that, without loss of generality, one can assume that each measure $\mu_n$ satisfies
\bel{wklmt1} Supp(\mu_n)~\subseteq~\Big\{ x\in \ov\Omega\,;~~
u_n(x)~\geq~c\alpha\,\kappa_0^{\alpha-1}\,|x|\Big\}~\doteq~A_n\,,
\eeq 	
where $c$ is the constant in (\ref{OP2}). 
Indeed, consider the decomposition
 $$\mu_n~=~\Hat \mu_n+\mu_n^*,$$
 where $\Hat\mu_n\doteq \chi_{A_n}\cdot \mu_n$ is concentrated on $A_n$, while $\mu^*_n$ is concentrated on
 $\ov\Omega\setminus A_n$. We notice that $A_n$ is a closed set, because
 $u_n$ is upper semicontinuous.
 

Observing that $u_n$ is a subsolution to the problem 
 \bel{Duf}\Delta u + f(x,u) -u\,\Hat  \mu_n~=~0\eeq
 with Neumann boundary conditions, we conclude that (\ref{Duf}), (\ref{NBC}) has a solution
$\Hat u_n \geq u_n$. For this solution, one has
 \bel{Hdif}\H(\Hat u_n,\Hat \mu_n)~\geq~\H(u_n,\mu_n)- \int_{\ov\Omega} u_n\,d\mu^*
 ~\geq~\H(u_n,\mu_n)- c\alpha\kappa_0^{\alpha-1}\int_{\ov\Omega} |x|\,d\mu^*.\eeq
Next, define the constants
$$\Hat \kappa_n ~= ~\Hat \mu_n(\ov\Omega),
\qquad \kappa^*_n ~=~ \mu_n^*(\ov \Omega),\qquad \kappa_n ~=~ \mu_n(\ov\Omega)~=~\hat\kappa_n + \kappa^*_n\,, $$
and consider an optimal irrigation plan for $\mu_n$, say
$\chi_n:[0, \kappa_n]\times [0,1]\mapsto\ov\Omega$.
By possibly performing a measure-preserving transformation of the domain $\Theta_n=[0, \kappa_n]$ into itself,
we can assume that the maps
$$\left\{\bega{l} \Hat \chi_n:[0, \Hat \kappa_n]\times [0,1]\mapsto\ov\Omega,\\[3mm]
\chi^*_n:[0, \kappa^*_n]\times [0,1]\mapsto\ov\Omega,\enda\right.
\qquad\qquad \left\{\bega{l} \Hat \chi_n(\theta,t)~=~ \chi_n(\theta,t),\\[3mm]
\chi_n^*(\theta,t) \doteq~\chi_n(\Hat\kappa_n+\theta,\, t),\enda\right.
$$
are admissible irrigation plans for $\Hat\mu_n$ and $\mu_n^*$, respectively (possibly not optimal).
We now have
\bel{Idif}\bega{l}\ds\I^\alpha(\mu_n)-\I^\alpha(\Hat \mu_n)~\geq~\int |x|_{\chi_n}^\alpha ~d\mathcal{H}^1-\int
|x|_{\strut\Hat \chi_n}^\alpha ~d\mathcal{H}^1~ \geq~\alpha\,\kappa_0^{\alpha-1}\int\Big( |x|_{\chi_n}-
|x|_{\Hat \chi_n}\Big) d\H^1(x)\\[4mm]
\ds
\qquad~=~\alpha\,\kappa_0^{\alpha-1}\int |x|_{\chi_n^*} d\H^1(x)
~\geq~\alpha\,\kappa_0^{\alpha-1}\int |x|\, d\mu_n^*.\enda
\eeq
Together, (\ref{Hdif}) and (\ref{Idif}) imply
$$\H(\Hat u_n,\Hat \mu_n)-c\I^\alpha(\Hat \mu_n)~\geq~\H(u_n,\mu_n)-c\I^\alpha(\mu_n).$$
By replacing each 
 pair $(u_n,\mu_n)$ with  $(\Hat u_n,\Hat \mu_n)$, we thus obtain a new maximizing sequence 
 for which (\ref{wklmt1}) holds.
\v
{\bf 3.} Using (\ref{wklmt1}), we now show that 
\bel{ubig}  Supp(\mu)~\subseteq~\Big\{ x\in \ov\Omega\,;~~
u(x)~\geq~c\alpha\,\kappa_0^{\alpha-1}\,|x|\Big\}.
\eeq 	
Indeed, assume that, on the contrary, 
there is a point $x_0\in Supp(\mu)$ such that 
\bel{choice}
u(x_0)~\leq~c\alpha \kappa^{\alpha-1}|x_0|-4\ve,\qquad \hbox{for some }      \ve>0.\eeq
By (\ref{limp}) there exists $t>0$ such that 
\bel{ch1}
u^{(t)}(x_0)+Kt~\leq~c\alpha \kappa^{\alpha-1}|x_0|-3\ve.\eeq
with $u^{(t)}$ defined as in (\ref{avp}). The continuity of $u^{(t)}$ implies
\bel{ch2}
u^{(t)}(x)+Kt~\leq~c\alpha \kappa^{\alpha-1}|x|-2\ve\eeq
for all $x\in B(x_0,r)\cap \ov\Omega$, with $r>0$ sufficiently small.
In turn, by the convergence $u^{(t)}_n(x)\to u^{(t)}(x)$, uniformly for all  $x\in \ov \Omega$, we have
\bel{ch3}
u_n^{(t)}(x)+Kt~\leq~c\alpha \kappa^{\alpha-1}|x|-\ve\eeq
for all $n\geq N_0 $ large enough and for all $x\in B(x_0,r)\cap \ov\Omega$.

By (\ref{wklmt1}), this implies that 
$$Supp(\mu_n)\cap B(x_0,r)~=~\emptyset,\qquad\forall n\geq N_0\,.$$
{}From the weak convergence $\mu_n\wto\mu$ it follows that $Supp(\mu)\cap B(x_0,r)~=~\emptyset$ as well,
contradicting the assumption $x_0\in Supp(\mu)$.
\v
{\bf 4.} Thanks to (\ref{ubig}) we can now define
 $$
 \mu^*~\doteq ~\frac{\nu}{u}.$$
By (\ref{wklmt}), 
 $u$ satisfies
\bel{u*eq}
\Delta u+f(x,u)-u\,\mu^*~=~0,\eeq
with Neumann boundary conditions (\ref{NBC}). 
 Following \cite{BCS}, we now establish the key inequality
\begin{equation}
\label{mu<}
\mu^*~\le~\mu.
\end{equation}
To prove that (\ref{mu<}) holds, thanks to the upper semicontinuity of $u$ it 
suffices to show that
\begin{equation}
\label{eq:inmu-1}
\int_{\overline\Omega}{\phi\over  \psi}\, d\nu~\le~
\int_{\overline\Omega}\phi\, d\mu,\qquad\hbox{for every }~\phi,\psi\in {\mathcal C}(\overline\Omega),~~\phi\ge0\,,~~~\psi >u.
\end{equation}

Since $\psi$ is continuous on the compact set $\ov\Omega$, 
we can choose $t,\delta>0$ small enough so that 
\bel{70}
u(x)~\leq ~u^{(t)}(x)+Kt ~<~\psi(x)
-\delta\qquad\qquad\forall x\in \ov\Omega\,.\eeq
By    (\ref{Gav}), as $n\to\infty$  the corresponding functions
$u_n^{(t)}$ converge to $u^{(t)}$ uniformly on $\ov\Omega$.
Hence for all $n$ large enough we have
$$u_n(x)~\leq~u_n^{(t)}(x)+Kt~\leq~u^{(t)}(x)+Kt +\delta~<~\psi(x)\qquad\qquad
\forall x\in \ov\Omega\,.$$
This yields
$$\int_{\ov\Omega}{\phi\over\psi}\, d\nu~=~\lim_{n\to\infty}
\int_{\ov\Omega}\phi\,{1\over\psi}\, d\nu_n~=~\lim_{n\to\infty}
\int_{\ov\Omega}\phi\,{u_n\over\psi}\, d\mu_n~\leq~\lim_{n\to\infty}
\int_{\ov\Omega}\phi\, d\mu_n~=~\int_{\ov\Omega}\phi\, d\mu\,,$$
proving (\ref{mu<}).

\v
{\bf 5.}  We conclude by proving the pair $(u,\mu^*)$ is optimal. Since $\{(u_n,\mu_n)\}_{n\geq1}$ is a maximizing sequence, using (\ref{wklmt}) and the lower semicontinuity of the irrigation cost $\I^\alpha$, one obtains
$$\bega{rl}
\ds \ov M &=\ds ~\lim_{n\to\infty}~
\Big[\H(u_n,\mu_n)-c\I^\alpha(\mu_n)\Big]\\[4mm]
\qquad &\leq\ds ~\lim_{n\to\infty}\int_{\Omega}f(x,u_n)~dx-c\liminf_{n\to\infty}\, \I^\alpha(\mu_n)~
\leq~\int_{\Omega}f(x,u)~dx-c\I^\alpha(\mu)\\[4mm]
\qquad &\leq~\ds \int_{\Omega}f(x,u)~dx-c\I^\alpha(\mu^{*}).
\enda$$
The last inequality follows from (\ref{mu<}) and
 the monotonicity of $\I^\alpha$, proved at (\ref{Ia2}). 
 By (\ref{mu<}) the weak convergence $\mu_n\wto\mu$ it follows
 $$\mu^*(\ov\Omega)~\leq ~\mu(\ov\Omega)~=~
 \kappa_0.$$ This completes the proof of the optimality of $(u,\mu^{*})$. 
\endproof
\v
We now prove an analogous existence result in the case of Dirichlet boundary conditions.

\v
\begin{thm}\label{thm5}  Let the assumptions (A1)-(A2) hold.
Then the 
maximization problem (\ref{OP2}), over all couples $(u,\mu)$ which satisfy (\ref{emu}), (\ref{DBC}), and (\ref{muO}), 
has an optimal solution.
\end{thm}
\v
{\bf Proof.} {\bf 1.} Call $\A$ the set of all admissible couples $(u,\mu)$ which satisfy (\ref{emu}), (\ref{DBC}), and (\ref{muO}). As in the previous case, the supremum $\ov M$ of the functional (\ref{OP2})
over all admissible couples $(u,\mu)\in\A$ satisfies (\ref{ovM}).
 Let $\{(u_n,\mu_n)\}_{n\geq 1}$ be a maximizing sequence. It is clearly not restrictive to 
 assume that $\mu_n\in \M_0$
 for every $n$.
\v
{\bf 2.} Let $w^*:\ov \Omega\mapsto [0,M]$ be the largest solution to the
elliptic problem with smooth coefficients
\bel{unc}\left\{\bega{rll}
\Delta w+ f(x,w)&=~0\qquad &x\in\Omega,\\[3mm]
w&=~0\qquad &x\in\partial\Omega.\enda\right.\eeq
By classical theory, $w^*$ can be constructed as the supremum of all functions $w:\ov\Omega\mapsto [0,M]$
which are subsolutions to (\ref{unc}).   Hence $w^*$ is well defined.

\begin{figure}[ht]
\centerline{\hbox{\includegraphics[width=7cm]{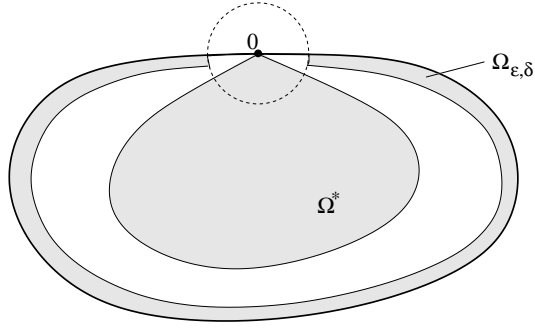}}}
\caption{\small The sets $\Omega^*$ and $\Omega_{\ve,\delta}$ considered in the proof of 
Theorem~\ref{thm5}.}
\label{f:ir13}
\end{figure}

For each $n\geq 1$, since $u_n\mu_n\geq 0$, by  Lemma \ref{l:comp}, the solution $u_n$ of (\ref{emu}),
(\ref{DBC}) satisfies
\bel{cuun}
u_n(x)~\leq ~w^*(x)\qquad\qquad \forall x\in \ov\Omega.\eeq
\v
{\bf 3.} Consider the set (see Fig.~\ref{f:ir13})
$$ \Omega^*~\doteq~\left\{x\in \Omega\,;~~w^*(x)\geq  c\,\alpha\kappa_0^{\alpha-1} |x|\right\}.
$$
Note that $\Omega^*$ is closed and
$$\Omega^*\cap\partial\Omega ~= ~\{0\}.$$
Denote by $\chi_{\Omega^*}$ the characteristic function of $\Omega^*$ and,
for each $n\geq 1$, consider the measure $\mu_n^*= \chi_{\Omega^*}\cdot\mu_n$ 
supported on $\Omega^*$.
Since $\mu_n^*\leq \mu_n$, by the comparison argument in Lemma \ref{l:comp}, 
we can find a solution $u_n^*$ of 
$$\Delta u+f(x, u) - u\,\mu^*_n~=~0$$
with Dirichlet boundary conditions (\ref{DBC}), such that 
 $$u_n~\leq ~u_n^*~\leq ~w^*.$$  We claim that $(u_n^*,\mu_n^*)_{n\geq 1}$ is another maximizing sequence. 
Indeed, 
\bel{we2g}
\H(u_n,\mu_n)-\H(u_n^*,\mu_n^*)~\leq~\int_{\Omega\setminus \Omega^*}
 u_n(x)\, d\mu_n~\leq~\int_{\Omega\setminus \Omega^*}
 c\,\alpha\kappa_0^{\alpha-1}|x| d\mu_n\,.
\eeq 
On the other hand, the same argument used at (\ref{Idif}) shows that the difference in the irrigation costs
can be estimated by
\bel{Idif2}\I^\alpha(\mu_n)-\I^\alpha(\mu_n^*)~\geq~\alpha \kappa_0^{\alpha-1}\int_{\Omega\setminus\Omega^*}
|x|\, d\mu_n\,.\eeq
Together, (\ref{we2g}) and (\ref{Idif2}) yield
\bel{final2}
\H(u_n,\mu_n)~-~\H(u_n^*,\mu_n^*)~\leq~c\I^\alpha(\mu_n)~-~c\I^\alpha(\mu_n^*),
\eeq
 proving that $\{(u_n^*,\mu_n^*)\}_{n\geq 1}$ is also a maximizing sequence.
 Without loss of generality, from now on we shall thus assume that 
 \bel{Spm}
 Supp(\mu_n)~\subseteq~\Omega^*\qquad\qquad\forall n\geq 1\,.\eeq
\v
{\bf 4.}  Consider
the sequence of measures $\nu_n\doteq u_n\,\mu_n$.
By possibly taking a subsequence, we can again assume that
(\ref{wklmt}) holds, for suitable positive measures $\mu,\nu$, supported  on $\Omega^*$. 
 Moreover,  for every fixed radius $r>0$, we can assume the convergence of the averaged values
\bel{avc} u_n^{(r)}(x)~\doteq~\avint_{B(x,r)\cap \Omega } u_n(y)\, dy~\to~\avint_{B(x,r)\cap \Omega } u(y)\, dy
~\doteq~u^{(r)}(x)\eeq
as $n\to\infty$, uniformly for $x\in \ov\Omega$.
\v
{\bf 5.}  We claim that 
\bel{ubig2} u(x)~\geq~c\alpha \kappa_0^{\alpha-1} |x|\qquad\forall x\in Supp(\mu).\eeq
Indeed, assume that, on the contrary, 
there is a point $x_0\in Supp(\mu)\subseteq \Omega^*$ such that 
\bel{ct0}
u(x_0)~\leq~c\alpha \kappa^{\alpha-1}|x_0|-4\ve,\qquad \hbox{for some }      \ve>0.\eeq
Clearly, this can hold only if $x_0\not= 0$.  Hence we can choose $r_0>0$ so that 
$B(x_0, 2r_0)\subset\Omega$.   

Since  $|f(x,u)|\leq K$,  all functions $u_n+{K\over 2d}|x|^2$, and $u+{K\over 2d}|x|^2$ are subharmonic on 
the open set $\Omega$.
Hence, there exists a constant $C$ such that, for every  $x\in B(x_0, r_0)$ and
$0<r\leq r_0$, all maps
$$r~\mapsto~ u^{(r)}_n(x) +Cr,\qquad\qquad r~\mapsto ~u^{(r)}_n(x) +Cr,$$ are nondecreasing.
Taking a sequence $r_k\to 0$, the pointwise values of $u_n, u$ can thus be defined as the infimum of a decreasing sequence of 
continuous functions:
\bel{uinf}u_n(x)~\doteq~\inf_{r>0} \avint_{B(x,r)} u_n(y)\, dy\,,\qquad\qquad 
u(x)~\doteq~\inf_{r>0} \avint_{B(x,r)} u(y)\, dy.\eeq
A contradiction is now achieved by the same argument used at (\ref{ch1})--(\ref{ch3}),
replacing the weigthed averages $u_n^{(t)}$ defined at (\ref{Gav}) with the standard averages $u_n^{(r)}$
in (\ref{avc}). 


\v

{\bf 6.}
By the previous step, 
we can define a measure $\mu^*$ supported on the open set $\Omega$, by setting
$$\mu^*~\doteq~{\nu\over u}\,.$$
Notice that, in principle, $\mu$ may contain a point mass at the origin.   In this case, to remove any ambiguity 
we define
$\mu^*(\{0\})=0$.    
By (\ref{wklmt}), the limit function
 $u$ satisfies (\ref{u*eq})
with Dirichlet boundary conditions (\ref{DBC}). 

The same arguments used in step {\bf 4} of the proof of Theorem~\ref{thm2} now show that $\mu^*\leq \mu$.
Hence the couple $(u, \mu^*)$ is admissible.   
\v
{\bf 7.} Finally, we prove that $(u, \mu^*)$ is optimal.  Indeed, 
on the set $\Omega\setminus\Omega^*$ all functions $u_n, u$ provide solutions to the 
semilinear elliptic equation with smooth coefficients
$$\Delta u + f(x,u)~=~0.$$
For any $\ve>0$, using  the Schauder regularity estimates \cite{E, F} up to the boundary, 
we can find $\rho>0$ 
such that all solutions $u_n$ are uniformly smooth on the set
$$\Omega_{\ve,\delta}~\doteq~\Big\{ x\in\Omega\,;~~dist(x,\partial \Omega)<\rho,~~|x|>\ve\Big\},$$
shown in Fig.~\ref{f:ir13}.
Hence, by Ascoli's theorem, by possibly taking a further subsequence we achieve the convergence
of the normal derivatives along the boundary
$$\partial_{\bfn(x)} u_n(x)~\to~\partial_{\bfn(x)} u(x)\qquad\qquad\forall  x\in\partial\Omega\setminus\{0\}.$$
Notice that, for any $\ve>0$, the convergence is uniform on the set
 $\partial\Omega\setminus
B(0,\ve)$.  
Observing that 
$0\leq u_n\leq w^*$ and similarly $0\leq u\leq w^*$, we deduce
\bel{pw}\partial_{\bfn(x)} u_n(x)~\leq~0,\qquad\qquad \partial_{\bfn(x)} w^*(x)~\leq~\partial_{\bfn(x)} u(x)~\leq~ 0.\eeq
Using (\ref{pw}), for any fixed $\ve>0$  one obtains
$$\bega{l}\ds\lim_{n\to\infty} \int_{\partial\Omega}\partial_\bfn u_n\, d\sigma
~\leq~\lim_{n\to\infty} \int_{\partial\Omega\setminus B(0,\ve)}\partial_\bfn u_n\, d\sigma
~=~\int_{\partial\Omega\setminus B(0,\ve)}\partial_\bfn u\, d\sigma\\[4mm]
\ds
~\leq~\int_{\partial\Omega}\partial_\bfn u\, d\sigma
-\int_{\partial\Omega\cap B(0,\ve)}\partial_\bfn w^*\, d\sigma.\enda
$$
Using the lower semicontinuity of the irrigation functional and the fact that $\mu^*\leq \mu$, we thus conclude
\bel{umopt}\bega{rl} \ov M&\ds =~\lim_{n\to\infty} \H(u_n,\mu_n)-\lim_{n\to\infty}c\I^\alpha(\mu_n)\\[4mm]
&\ds\leq~\lim_{n\to\infty} \int_\Omega f(x, u_n)\, dx +\lim_{n\to\infty} \int_{\partial\Omega}\partial_\bfn u_n\, d\sigma
-c\I^\alpha(\mu)\\[4mm]
&\ds\leq~ \int_\Omega f(x, u)\, dx +\int_{\partial\Omega}\partial_\bfn u\, d\sigma
+\int_{\partial\Omega\cap B(0,\ve)}\partial_\bfn w^*\, d\sigma
-c\I^\alpha(\mu^*)\\[4mm]
&\ds=~\H(u,\mu^*) -c\I^\alpha(\mu^*) -\int_{\partial\Omega\cap B(0,\ve)}\partial_\bfn w^*\, d\sigma.\enda
\eeq
By choosing $\ve>0$ small, the last integral on the right hand side of (\ref{umopt}) can be made 
arbitrarily small.   Hence $\H(u,\mu^*) -c\I^\alpha(\mu^*)~\geq~\ov M$, proving the optimality of $(u,\mu^*)$.
\endproof

\v
\section{Concluding remarks}

In this paper we assumed that the primary goal of 
tree leaves (tree roots) is to gather sunlight 
(water and nutrients from the soil, respectively).
We then tried to determine shapes that most efficiently achieve these goals.
The search for these optimal shapes has been formulated
as a maximization problem for certain functionals, in the spirit of the classical Calculus of Variations \cite{ABM}. 

While our present analysis is purely theoretical, optimal shapes 
may be computed by the numerical algorithms recently 
developed in \cite{M, OS, PSX, S2}.
It will then be of interest to compare numerical simulations with 
the shapes actually observed in nature.    In this direction, we expect that
root shapes which maximize our harvest functional will look very similar to the 
actual roots of
biological trees.     

On the other hand, we guess that in some cases
the shapes which maximize the gathered sunlight will resemble an optimal disposition
of solar panels, more than actual tree branches.
If this is the case, it would indicate that the efficiency in capturing sunlight 
has not been the primary goal driving the evolution of  plant shapes.
In computer simulations of tree growth \cite{AK, H, RLP}, 
the most realistic images are produced
by algorithms based on the idea of conquering space. 
This suggests that tree shapes have evolved as the result of
a competitive game among plants, rather than an optimization problem.
A mathematical modeling of such a game remains to be worked out.

\v
{\bf Acknowledgment.}
This research was partially supported by NSF with  
grant DMS-1714237, ``Models of controlled biological growth".

\v


\begin{thebibliography}{6111}

\bibitem{AFP}
L.~Ambrosio, N.~Fusco, and D.~Pallara,
{\it Functions of Bounded Variation and Free Discontinuity Problems.}
Oxford University Press, 2000.

\bibitem{AGS} L.~Ambrosio, N.~Gigli, and G.~Savar\'e, 
{\it Gradient Flows in Metric Spaces and in the Space of Probability Measures.} 
Birkh\"auser, Basel, 2005.


\bibitem{AK} M~Aono and T~L.~Kunii, Botanical tree image generation. 
{\it IEEE Computer Graphics and Appl.} {\bf 4, 5} (1984), 10--34.

\bibitem{AG}
D.~H.~Armitage and S.~J.~Gardiner,
{\it  Classical Potential Theory.}
 Springer-Verlag, London, 2001.


\bibitem{ABM} H.~Attouch, G.~Buttazzo and G.~Michaille,
{\it Variational Analysis in Sobolev and BV Spaces: Applications to PDEs and Optimization}, Second Edition.
MOS-SIAM Series on Optimization, 2014.


\bibitem{BM}~J.~M.~Ball and F.~Murat, Remarks on Chacon's biting lemma.
{\it Proc. Amer. Math. Soc.} {\bf 107} (1989),
655--663.

\bibitem{BCM} M.~Bernot,  V.~Caselles, and J.~M.~Morel,
{\it  Optimal transportation networks. Models and theory.} Springer Lecture Notes in Mathematics {\bf 1955},
Berlin, 2009.

\bibitem{Bi}
P.~Billingsley, {\it Convergence of Probability Measures.} Wiley, New York, 1999.
 
\bibitem{BG1}
 L.~Boccardo and  T.~Gallou\"et,
 Non-linear elliptic and parabolic equations involving measure data.
{\it J.~Functional Analysis} {\bf  87} (1989), 149--169.
 
 \bibitem{BGO}
 L.~Boccardo and  T.~Gallou\"et, and L.~Orsina,
 Existence and uniqueness of entropy solutions for nonlinear elliptic equations with measure data.
{\it Ann. Institut H.~Poincar\'e Nonlin. Anal.} {\bf 13} (1996), 539--551.
 
 
\bibitem{BS} L.~Brasco and F.~Santambrogio, 
An equivalent path functional formulation of branched transportation problems. 
{\it Discrete Contin. Dyn. Syst.} {\bf 29} (2011), 845--871.

\bibitem{BCS} A.~Bressan,
G.~Coclite and W.~Shen,
A multi-dimensional optimal harvesting problem with measure valued solutions,
{\it SIAM J. Control Optim.} {\bf 51} (2013), 1186--1202.


  
  

\bibitem{BC80}
J.~K.~Brooks and R.~V.~Chacon,
Continuity and compactness of measures. 
{\it Adv. in Math.} {\bf  37} (1980), 16--26. 

\bibitem{DMOP} G.~Dal Maso, F.~Murat, L.~Orsina, and A.~Prignet,
Renormalized solutions of elliptic equations with general measure data
{\it Ann. Scuola Norm. Sup. Pisa Cl. Sci.}  {\bf 28} (1999), 741--808.

\bibitem{E}
L.~C.~Evans,  {\it   Partial Differential Equations.  Second edition. }
American Mathematical Society, Providence, RI, 2010.

\bibitem{EG}
L.~C.~Evans and Ronald~F.~Gariepy,  {\it  Measure Theory and Fine Properties of Functions.  }
CRC Press, 1991. 

\bibitem{FZ} H.~Federer and W.~Ziemer,
The Lebesgue set of a function whose distribution derivatives are p-th power summable.
{\it  Indiana Univ. Math. J.} {\bf 22} (1972), 139--158.

\bibitem{F} A.~Friedman, {Partial differential equations of parabolic type.} Prentice-Hall, 
Englewood Cliffs, N.J., 1964.

\bibitem{GMS}
M.~Giaquinta, G.~Modica, and J.~Soucek, 
{\it Cartesian Currents in the Calculus of Variations I}.  Springer-Verlag, Berlin, 1998.


\bibitem{G}
E.~N.~Gilbert. Minimum cost communication networks. 
{\it Bell System Tech. J.} {\bf 46} (1967), 2209--2227.


\bibitem{H} 
H.~Honda, Description of the form of trees by the parameters of the tree-like body
{\it J.~Theor. Biology} {\bf 31} (1971), 331--338.

\bibitem{MMS} F.~Maddalena, J.~M.~Morel, and S.~Solimini, 
A variational model of irrigation patterns, 
{\it Interfaces Free Bound.} {\bf 5}  (2003), 391--415.

\bibitem{MS} F.~Maddalena and S.~Solimini,
Synchronic and asynchronic descriptions of irrigation problems. 
{\it Adv. Nonlinear Stud.} {\bf 13} (2013), 583--623. 

\bibitem{M} A.~Monteil,
Uniform estimates for a Modica-Mortola type approximation of branched transportation, 
{\it ESAIM Control Optim. Calc. Var.} {\bf 23}  (2017), 309--335.



 \bibitem{OS}
 E.~Oudet and F.~Santambrogio,
 A Modica-Mortola approximation for branched transport and applications. 
{\it  Arch. Rational Mech.  Anal.} {\bf  201} (2011), 115--142.


\bibitem{PSX} P.~Pegon, F.~Santambrogio, and Q.~Xia, 
A fractal shape optimization problem in branched transport. Preprint 2017.

\bibitem{RLP}
A.~Runions, B.~Lane, and P.~Prusinkiewicz,
Modeling Trees with a Space Colonization Algorithm,
{\it Eurographics Workshop on Natural Phenomena}, 2007.


\bibitem{S2}  F.~Santambrogio,  A Modica-Mortola approximation for branched transport.   
{\it C. R. Acad. Sci. Paris, Ser. I},
{\bf  348} (2010) 941--945.

\bibitem{X3}
Q.~Xia, Optimal paths related to transport problems, 
{\it Comm. Contemp. Math.} {\bf 5}  (2003), 251--279.


\bibitem{X15} Q.~Xia, 
Motivations, ideas and applications of ramified optimal transportation.
{\it ESAIM, Math. Model. Numer. Anal.}
{\bf  49} (2015),  1791--1832. 





\end{thebibliography}
\end{document}